\documentclass[12pt, reqno]{amsart}
\usepackage{a4, amsmath, amsthm, amssymb}
\usepackage{mathbbol}
\usepackage{txfonts}
\usepackage{mathrsfs}

\newtheorem{thm}{Theorem}[section]
\newcommand{\bt}{\begin{thm}}
\newcommand{\et}{\end{thm}}

\newtheorem{lem}[thm]{Lemma}
\newcommand{\bl}{\begin{lem}}
\newcommand{\el}{\end{lem}}

\newtheorem{prop}[thm]{Proposition}
\newcommand{\bp}{\begin{prop}}
\newcommand{\ep}{\end{prop}}

\newtheorem{defn}[thm]{Definition}
\newcommand{\bd}{\begin{defn}}      
\newcommand{\ed}{\end{defn}}

\newcommand{\thmref}[1]{Theorem~\ref{#1}}
\newcommand{\secref}[1]{Section~\ref{#1}}
\newcommand{\lemref}[1]{Lemma~\ref{#1}}
\newcommand{\defref}[1]{Definition~\ref{#1}}
\newcommand{\propref}[1]{Proposition~\ref{#1}}
\newcommand{\N}{\mathbb{N}}
\newcommand{\R}{\mathbb{R}}
\newcommand{\Z}{\mathbb{Z}}
\newcommand{\id}{\operatorname{id}}

\newcommand{\hm}{{\mathscr H}} 
\newcommand{\lm}{{\mathscr L}} 

\newcommand{\clB}{B}       

\newcommand{\Lip}{\operatorname{Lip}}
\newcommand{\Lipb}{\Lip^\text{\rm b}}       
\newcommand{\Lipbs}{\Lip_\text{\rm bs}}     
\newcommand{\Liploc}{\Lip_\text{\rm loc}}   
\newcommand{\Bor}{{\mathscr B}}             
\newcommand{\Borb}{\Bor^\infty}             
\newcommand{\Borbbs}{\Borb_\text{\rm bs}}   
\newcommand{\Borbloc}{\Borb_\text{\rm loc}} 

\newcommand{\form}[2]{\Lipbs(#1)\times[\Liploc(#1)]^{#2}}

\newcommand{\LWc}[2]{{\mathbf M}_{\text{\rm loc},\,#2}(#1)}  
\newcommand{\LWnc}[2]{{\mathbf N}_{\text{\rm loc},\,#2}(#1)} 
\newcommand{\LWirc}[2]{{\mathscr I}_{\text{\rm loc},\,#2}(#1)}
\newcommand{\LWic}[2]{{\mathbf I}_{\text{\rm loc},\,#2}(#1)}

\newcommand{\mass}{{\mathbf M}}

\newcommand{\AKc}{{\mathbf M}}    
\newcommand{\AKnc}{{\mathbf N}}   
\newcommand{\AKirc}{{\mathscr I}} 
\newcommand{\AKic}{{\mathbf I}}   

\newcommand{\flatnorm}{{\mathscr F}}

\newcommand{\rstr}{\:\mbox{\rule{0.1ex}{1.2ex}\rule{1.1ex}{0.1ex}}\:}
\newcommand{\bdry}{\partial}
\newcommand{\on}[1]{|_{#1}}
\newcommand{\spt}{\operatorname{spt}}

\begin{document}

\title[Pointed flat convergence]{The pointed flat compactness theorem\\for locally integral currents}

\author{Urs Lang}

\address{Department of Mathematics\\
ETH Zurich\\
8092 Zurich, Switzerland}
\email{lang@math.ethz.ch}

\author{Stefan Wenger}

\address{Department of Mathematics\\
University of Illinois at Chicago\\
851 S. Morgan Street\\
Chicago, IL 60607--7045, USA}
\email{wenger@math.uic.edu}

\date{\today}

\thanks{The second author is partially supported by NSF grant DMS 0956374}


\begin{abstract}
Recently, a new embedding/compactness theorem for integral currents in a sequence of metric spaces has been established by the second author. We present a version of this result for locally integral currents in a sequence of pointed metric spaces. To this end we introduce another variant of the Ambrosio--Kirchheim theory of currents in metric spaces, including currents with finite mass in bounded sets.
\end{abstract}

\maketitle

\section{Introduction}

In the recent article~\cite{Wenger-cpt}, the second-named author proved a new compactness theorem that combines features of two other powerful results: Gromov's theorem for uniformly compact families of metric spaces and the compactness theorem for integral currents from geometric measure theory. When applied to a sequence $M_n$ of compact, connected and oriented Riemannian $m$-manifolds\footnote{In this paper, the term ``manifold'' allows for boundaries.}, the theorem says that if the diameters, the volumes and the volumes of the boundaries are uniformly bounded, then there exist a subsequence $M_{n_j}$, a complete metric space $Z$, and isometric (i.e.~distance preserving) embeddings $\varphi_j \colon M_{n_j} \hookrightarrow Z$ such that the images $\varphi_j(M_{n_j})$, viewed as integral currents $\varphi_{j\#}[M_{n_j}]$, converge in the flat topology to an integral current $T$ in $Z$. Here the terms ``integral current'' and ``flat topology'' are understood in the sense of the metric theory of currents introduced by Ambrosio--Kirchheim~\cite{AK-c}. In case the $M_n$ have no boundary, there is a stronger conclusion: As $j  \to \infty$, the filling volume of the cycles $T - \varphi_{j\#}[M_{n_j}]$ tends to zero, i.e.~they bound $(m+1)$-currents with smaller and smaller mass. The general formulation of the theorem refers to sequences of integral currents $T_n$ in complete metric spaces $X_n$. The purpose of the present article is to provide a ``pointed version'' of this result for locally integral currents, akin to Gromov's embedding/compactness theorem for pointed metric spaces \cite{Gromov-polynomialgrowth}.

The metric currents of~\cite{AK-c} have finite mass by definition and are therefore not suited for the envisaged pointed compactness theorem. In~\cite{Lang-c}, the first-named author presented a theory of general metric currents with locally compact supports, comprising currents $T$ with locally finite mass, whose measure $\|T\|$ is a Radon measure. However, this setup does again not fit our purpose, as the local compactness may fail to persist in the limit. We therefore present yet another variation of the theory of currents in a metric space $Z$. This will be discussed at some detail in \secref{Section:local-currents}; here we just briefly explain the notation needed for the main results. The abelian group of $m$-dimensional locally integral currents in $Z$ will be denoted by $\LWic{Z}{m}$. The measure $\|T\|$ associated with an element $T$ is finite on bounded sets and furthermore concentrated on some countably $m$-rectifiable set. We emphasize that throughout this article the subscript ``loc'' refers to a finiteness property on all bounded sets rather than on suitable point neighborhoods. The support in $Z$ and boundary of a current $T \in \LWic{Z}{m}$ will be denoted by $\spt T$ and $\bdry T$; the latter is an element of $\LWic{Z}{m-1}$. For instance, every connected and oriented Riemannian $m$-manifold $M$ that is complete as a metric space induces a current $[M] \in \LWic{M}{m}$ with $\spt[M] = M$ and $\bdry[M] = [\bdry M]$. The push-forward $\varphi_\#T$ of $T \in \LWic{Z}{m}$ is defined and belongs to $\LWic{Z'}{m}$ whenever $\varphi$ maps $\spt T$ into $Z'$ such that preimages of bounded sets are bounded and $\varphi$ is Lipschitz on bounded subsets of $\spt T$. The boundary of $\varphi_\#T$ equals $\varphi_\#(\bdry T)$. We say that a sequence $(T_j)$ in $\LWic{Z}{m}$ converges in the local flat topology to a current $T \in \LWic{Z}{m}$ if for every bounded closed set $B \subset Z$ there is a sequence $(S_j)$ in $\LWic{Z}{m+1}$ such that 
\[
(\|T - T_j - \bdry S_j\| + \|S_j\|)(B) \to 0;
\]
in other words, $T - T_j$ can be expressed as $R_j + \bdry S_j$, with $R_j \in \LWic{Z}{m}$ and $S_j \in \LWic{Z}{m+1}$, such that $(\|R_j\| + \|S_j\|)(B) \to 0$. Then $\bdry T_j \to \bdry T$ in the local flat topology of $\LWic{Z}{m-1}$.

Given a metric space $X$ and $x_0 \in X$, we denote by $\clB(x_0,r) := \left\{x \in X: d(x_0,x)\leq r\right\}$ the closed ball of radius $r$ at $x_0$.

The main result of this article is the following pointed version of~\cite[Theorem~1.2]{Wenger-cpt}. The proof uses the same decomposition techniques and will be given in \secref{subsection:ptd-cptness}.

\bt\label{thm:main}
Let $(X_n)$ be a sequence of complete metric spaces, $x_n \in X_n$, and let $T_n\in\LWic{X_n}{m}$, $m \ge 1$. Suppose that for every $r>0$,
\begin{equation*}
\sup_n\left[ \|T_n\|(\clB(x_n,r)) + \|\bdry T_n\|(\clB(x_n,r)) \right] < \infty.
\end{equation*}
Then there exist a subsequence $(n_j)$, a complete metric space $Z$ with basepoint $z_0$, and isometric embeddings $\varphi_j\colon X_{n_j} \hookrightarrow Z$ such that $\varphi_j(x_{n_j}) \to z_0$ and $(\varphi_{j\#}T_{n_j})$ converges in the local flat topology to some $T\in\LWic{Z}{m}$.
\et

Similarly as in~\cite{Wenger-cpt}, a somewhat stronger conclusion holds if $\bdry T_n = 0$ for all $n$, or, more generally, if for every $r>0$, $\spt(\bdry T_n) \cap \clB(x_n,r) = \emptyset$ for almost every $n$. Then, for every bounded set $B \subset Z$, $\spt(\bdry(\varphi_{j\#}T_{n_j})) \cap B = \emptyset$ for almost all $j$, in particular $\bdry T = 0$. In this situation, $T - \varphi_{j\#}T_{n_j}\to 0$ in the ``local filling sense'': For every bounded closed set $B \subset Z$ there is a sequence $(S'_j)$ in $\LWic{Z}{m+1}$ such that $\spt(T - \varphi_{j\#}T_{n_j} - \bdry S'_j) \cap B = \emptyset$ for almost all $j$, and 
\[
\|S'_j\|(B)  \to 0.
\]
This will be discussed in \secref{subsection:fill-conv}.

The next result shows that the limit given by \thmref{thm:main} is unique, up to a pointed isometry. 

\bp\label{prop:intro-uniqueness}
Let $(X_n)$ be a sequence of complete metric spaces, $x_n \in X_n$, and let $T_n\in\LWic{X_n}{m}$. Suppose there exist two complete metric spaces $Z$, $Z'$ with basepoints $z_0$, $z'_0$ and isometric embeddings $\varphi_n\colon X_n\hookrightarrow Z$, $\varphi'_n\colon X_n\hookrightarrow Z'$ such that $\varphi_n(x_n) \to z_0$, $\varphi'_n(x_n) \to z'_0$, $(\varphi_{n\#}T_n)$ converges in the local flat topology to $T\in\LWic{Z}{m}$, and $(\varphi'_{n\#}T_n)$ converges in the local flat topology to $T'\in\LWic{Z'}{m}$. Then there is an isometry $\psi \colon \{z_0\} \cup \spt T \to \{z'_0\} \cup \spt T'$ with $\psi(z_0) = z'_0$ and $\psi_\#T = T'$.
\ep

This will be proved in~\secref{subsection:uniqueness}. See~\cite[Theorem~1.3]{Wenger-cpt} for the analog in the bounded case.

Finally, in \secref{subsection:ultralimit}, we shall discuss the following proposition, relating the above results to other types of pointed limits. Compare~\cite[Proposition~2.2]{Wenger-cpt}. Here it suffices to assume that $\varphi_{n\#}T_n \to T$ weakly, i.e., pointwise as functionals. Convergence in the local flat topology implies weak convergence, and the reverse implication holds under mild additional assumptions, cf.~\cite{Wenger-flatconv}. 

\bp\label{prop:ultralimit}
For $n \in \N = \{1,2,\dots\}$, let $X_n$ be a complete metric space, $x_n \in X_n$, and let $T_n\in\LWic{X_n}{m}$. Suppose there exist a complete metric space $Z$ with basepoint $z_0$ and isometric embeddings $\varphi_n\colon X_n\hookrightarrow Z$ such that $\varphi_n(x_n) \to z_0$ and $(\varphi_{n\#}T_n)$ converges weakly to $T\in\LWic{Z}{m}$. 
\begin{enumerate}
\item
For every non-principal ultrafilter $\omega$ on $\N$, there is an isometric embedding of $\{z_0\} \cup \spt T$ into the ultralimit $(X_\omega,x_\omega) = \lim_\omega(X_n,x_n)$ that maps $z_0$ to $x_\omega$.
\item
If $(X_n,x_n)$ converges in the pointed Gromov--Hausdorff sense to a pointed proper metric space $(Y,y_0)$, then there is an isometric embedding of $\{z_0\} \cup \spt T$ into $Y$ that maps $z_0$ to $y_0$.
\end{enumerate}
\ep

The metric approach to currents, originally proposed by De Giorgi~\cite{DeGiorgi-funzionali}, employs $(m+1)$-tuples of real-valued functions as a substitute for differential $m$-forms. If the underlying metric space possesses a smooth structure, the $m$-form $f d\pi_1 \wedge \ldots \wedge d\pi_m$ is represented by the tuple $(f,\pi_1,\dots,\pi_m)$. In the theory of currents of finite mass presented in~\cite{AK-c}, for a complete metric space $Z$, the tuples consist of Lipschitz functions on $Z$, and the first entry is bounded in addition. An $m$-dimensional current is then an $(m+1)$-linear functional on $\Lipb(Z) \times [\Lip(Z)]^m$ satisfying a set of further conditions, depending on the class of currents under consideration. In the theory of local metric currents exposed in~\cite{Lang-c}, the underlying metric space is at first assumed to be locally compact, and the first function of a test tuple is Lipschitz with compact support, the remaining ones are locally Lipschitz. A natural extension of the theory then includes currents with locally compact support in an arbitrary metric space $Z$. The ``boundedly finite'' theory of metric currents discussed here uses ``forms'' in $\form{Z}{m}$, where ``bs'' stands for ``bounded support'' and $\Liploc(Z)$ denotes the space of functions that are Lipschitz on bounded sets. We point out again that here the subscript ``loc'' has a different meaning than in~\cite{Lang-c} unless $Z$ is proper, i.e.~bounded closed subsets of $Z$ are compact. We shall discuss the fundamentals of the theory in detail, so that no prior knowledge of~\cite{AK-c} or~\cite{Lang-c} is required. However, some of the more profound results, such as the boundary rectifiability and closure theorems, will be deduced from their analogues in~\cite{AK-c}.


\section{Metric currents with finite mass on bounded sets}\label{Section:local-currents}

Let $Z$ and $Z'$ be metric spaces. We denote by $\Lip(Z,Z')$ the set of all Lipschitz maps from $Z$ to $Z'$ and by $\Liploc(Z, Z')$ the set of all maps from $Z$ to $Z'$ that are Lipschitz continuous on bounded subsets of $Z$. We write $\Lip(Z)$ and $\Liploc(Z)$ for the vector spaces $\Lip(Z,\R)$ and $\Liploc(Z,\R)$, respectively. Note that the latter is an algebra. Furthermore, $\Lipb(Z)$ denotes the algebra of bounded real-valued Lipschitz functions on $Z$ and $\Lipbs(Z)$ the subalgebra of functions with bounded support. The Lipschitz constant of a function $f$ is denoted by $\Lip(f)$.

\subsection{Metric functionals}

We first consider real-valued functions on the space of $(m+1)$-tuples $\form{Z}{m}$, where $m \ge 0$. A typical such tuple will be denoted by $(f,\pi_1,\dots,\pi_m)$, and we may use $(f,\pi)$ as a shorthand. In case $m = 0$, $\form{Z}{m}$ should be read as $\Lipbs(Z)$ and $(f,\pi)$ as $f$. Let now 
$$T \colon \form{Z}{m} \to \R$$ 
be given. For any tuple $(g,\tau) := (g,\tau_1,\dots,\tau_l) \in [\Liploc(Z)]^{l+1}$ with $0\leq l\leq m$, we define a function $T\rstr(g,\tau) \colon \form{Z}{m-l} \to \R$ by 
\begin{equation*}
 T\rstr(g,\tau)\: (f,\pi_1,\dots,\pi_{m-l}):= T(fg,\tau_1,\dots,\tau_l,\pi_1,\dots,\pi_{m-l})
\end{equation*}
and call it the restriction of $T$ to $(g,\tau)$. 
For any map $\varphi\in\Liploc(Z,Z')$ with the property that $\varphi^{-1}(A)$ is bounded for every bounded set $A\subset Z'$, we define a function $\varphi_{\#}T \colon \form{Z'}{m} \to \R$ by
\begin{equation*}
 \varphi_{\#}T\,(f,\pi_1,\dots,\pi_m):= T(f\circ\varphi,\pi_1\circ\varphi, \dots, \pi_m\circ\varphi)
\end{equation*}
and call it the push-forward under $\varphi$ of $T$. 

\bd\label{def:metric-fctl}
A function $T\colon\form{Z}{m}\to\R$, $m \ge 0$, is called an $m$-dimensional metric functional on $Z$ if the following properties hold:
 \begin{enumerate}
  \item $T$ is multilinear;
  \item $T$ is continuous in the following sense: If $\pi_i, \pi_i^j\in\Liploc(Z)$, $\pi_i^j\to\pi_i$ pointwise everywhere as $j\to\infty$ and $\sup_{i,j}\Lip(\pi_i^j\on{A})<\infty$ for every bounded set $A\subset Z$, then
   \begin{equation*}
    T(f, \pi_1^j,\dots, \pi_m^j)\to T(f, \pi_1, \dots,\pi_m);
   \end{equation*}
  \item $T$ is local in the following sense: If there exist $i$ and $\delta > 0$ such that $\pi_i$ is constant on $\left\{z: d(z,\spt f) \le \delta\right\}$, then $T(f,\pi_1,\dots,\pi_m)=0$.
 \end{enumerate}
\ed

A $0$-dimensional metric functional on $Z$ is just a linear functional on $\Lipbs(Z)$. It is not difficult to verify that restrictions and push-forwards of metric functionals are metric functionals. To check property~(iii) for $\varphi_\#T$, observe that since $\varphi$ is Lipschitz on tubular neighborhoods of the bounded set $\spt(f \circ \varphi)$, for every $\delta' > 0$ there is a $\delta > 0$ such that $\varphi$ maps $\left\{z: d(z,\spt(f \circ \varphi)) \le \delta\right\}$ into $\left\{z': d(z',\spt f) \le \delta'\right\}$.
As a consequence of the defining conditions of a metric functional, the locality property also holds in a strict form:

\bl\label{lem:strict-loc}
If some $\pi_i$ is constant on $\spt f$, then $T(f,\pi_1,\dots,\pi_m)=0$.
\el

\begin{proof}
Suppose first that $\pi_i\on{\spt f} = 0$ for some $i$. For $j \in \N$, define $\beta_j\colon\R \to \R$ so that $\beta_j(s) = \max\{0,s-j^{-1}\}$ for $s \ge 0$ and $\beta_j(s) = \min\{0,s+j^{-1}\}$ for $s \le 0$. Note that $\beta_j$ is $1$-Lipschitz. As $j \to \infty$, $\beta_j \circ \pi_i$ converges pointwise to $\pi_i$. It thus follows from the continuity property of $T$ that 
 \begin{equation*}
  T(f,\pi_1,\dots,\pi_m) = \lim_{j\to\infty} T(f,\pi_1,\dots,\pi_{i-1},\beta_j \circ \pi_i,\pi_{i+1},\dots,\pi_m).
 \end{equation*}
Furthermore, since $\pi_i\on{\spt f} = 0$ and $\pi_i$ is Lipschitz on tubular neighborhoods of the bounded set $\spt f$, for every $j$ there is a $\delta_j > 0$ such that $|\pi_i(z)| \le j^{-1}$ whenever $d(z,\spt f) \le \delta_j$. Then $(\beta_j \circ \pi_i)(z) = 0$ for all such $z$, thus the above equality and the locality of $T$ imply $T(f,\pi_1,\dots,\pi_m) = 0$. 

To conclude the proof in the general case, note that by~(i) and~(iii), adding a constant to one of $\pi_1,\dots,\pi_m$ does not change the value $T(f,\pi_1,\dots,\pi_m)$. 
\end{proof}

Now let $T$ be a metric functional of dimension $m \ge 1$ on $Z$. We define its boundary $\bdry T \colon \form{Z}{m-1} \to \R$ by
\begin{equation*}
 \bdry T(f,\pi_1,\dots,\pi_{m-1}):= T(\sigma,f,\pi_1,\dots,\pi_{m-1}),
\end{equation*}
where $\sigma\in\Lipbs(Z)$ is any function satisfying $\sigma\on{\spt f} = 1$. If $\sigma'$ is another such function, then $f$ vanishes on $\{\sigma \ne \sigma'\}$ and hence on $\spt(\sigma - \sigma')$, so $T(\sigma -\sigma',f,\pi_1,\dots,\pi_{m-1}) = 0$ by the above lemma. Hence $\bdry T$ is well-defined. Clearly $\bdry T$ satisfies properties~(i) and~(ii) of \defref{def:metric-fctl}. To verify~(iii), suppose that for some $i \in \{1,\dots,m-1\}$, $\pi_i$ is constant on a tubular neighborhood of $\spt f$. Then it is possible to choose $\sigma \in \Lipbs(Z)$ with $\sigma\on{\spt f} = 1$ such that $\pi_i$ is constant on some tubular neighborhood of $\spt\sigma$ and hence $\bdry T(f,\pi_1,\dots,\pi_{m-1}) = T(\sigma,f,\pi_1,\dots,\pi_{m-1}) = 0$ by the locality of $T$. Thus $\bdry T$ is an $(m-1)$-dimensional metric functional. It is easy to check that 
\begin{equation*}
 \varphi_\#(\bdry T) = \bdry (\varphi_\# T).
\end{equation*}
Furthermore, if $m \ge 2$, then
\begin{equation*}
 \bdry(\bdry T) = 0.
\end{equation*}
To see this, let $(f,\pi) := (f,\pi_1,\dots,\pi_{m-2}) \in \form{Z}{m-2}$ and choose $\varrho,\sigma,\tau \in \Lipbs(Z)$ such that $\varrho\on{\spt f} = 1$, $\sigma\on{\spt\varrho} = 1$, and $\tau\on{\spt\sigma} = 1$, in particular $\sigma\on{\spt f} = 1$. By definition,
\begin{equation*}
 \bdry(\bdry T)(f,\pi) = \bdry T(\sigma,f,\pi) = T(\tau,\sigma,f,\pi).
\end{equation*}
Now $f$ vanishes on $\{\tau \ne \varrho\}$ and hence on $\spt(\tau - \varrho)$. It follows that $T(\tau,\sigma,f,\pi) = T(\varrho,\sigma,f,\pi)$, and this last term is zero since $\sigma\on{\spt\varrho} = 1$.

\bp
Every metric functional of dimension $m \ge 2$ on $Z$ is alternating in the $m$ arguments $\pi_1,\dots,\pi_m \in \Liploc(Z)$.
\ep
 
\begin{proof}
This is shown by the same argument as in the proofs of \cite[Theorem~3.5]{AK-c} and \cite[Proposition~2.4]{Lang-c}.
\end{proof}

\subsection{Mass}

We now introduce the mass of a multilinear functional and then discuss metric functionals with finite mass on bounded sets. The local mass bound implies a stronger continuity property that involves the first argument of the functional. This leads to further properties, justifying the use of the term ``current''. In~\cite{AK-c}, the measure associated with a current of finite mass is concentrated on a $\sigma$-compact set. We bring this property into play at an early stage (cf.~\propref{prop:sigma-cpt}), as a substitute for the local compactness underlying the approach of~\cite{Lang-c}.

We denote by $\Lip_1(Z) \subset \Lip(Z)$ the set of all $1$-Lipschitz functions. For every multilinear function $T \colon \form{Z}{m} \to \R$ and every open set $V \subset Z$ we define the mass of $T$ in $V$ as the possibly infinite quantity
\[
 \mass_V(T) := \sup \sum_{\lambda \in \Lambda} T(f_\lambda,\pi_\lambda),
\]
where the supremum is taken over all finite families $\bigl((f_\lambda,\pi_\lambda)\bigr)_{\lambda \in \Lambda}$ such that $(f_\lambda,\pi_\lambda) = (f_\lambda,\pi_{\lambda,1},\dots,\pi_{\lambda,m}) \in \Lipbs(Z) \times [\Lip_1(Z)]^m$, $\spt f_\lambda \subset V$, and $\sum_{\lambda \in \Lambda}|f_\lambda| \le 1$. In case $m = 0$, 
\[
 \mass_V(T) = \sup\left\{T(f): f \in \Lipbs(Z),\,\spt f \subset V,\,|f| \le 1\right\}.
\]
If a sequence $(T_n)$ of multilinear functions converges pointwise on $\form{Z}{m}$ to a multilinear function $T$, then clearly
\[
 \mass_V(T) \le \liminf_{n \to \infty} \mass_V(T_n)
\]
for every open set $V \subset Z$, i.e.~$\mass_V$ is lower semicontinuous. Pointwise convergence of multilinear functions will be referred to as weak convergence. We write $\mass(T) := \mass_Z(T)$ for the total mass. 
We now define a set function $\|T\| \colon 2^Z \to [0,\infty]$ by
\[
 \|T\|(A) := \inf\left\{\mass_V(T): \text{$V \subset Z$ is open, $A \subset V$}\right\}.
\]
If $A$ is open, then obviously $\|T\|(A) = \mass_A(T)$. For two multilinear functions $T,T' \colon \form{Z}{m} \to \R$ and $r \in \R$ we have 
\[
\|T +T'\| \le \|T\| + \|T'\|, \quad \|r T\| = |r| \,\|T\|. 
\]
Under a suitable $\sigma$-compactness assumption, $\|T\|$ turns out to be an outer measure.

\bp\label{prop:sigma-cpt}
Let $T \colon \form{Z}{m} \to \R$ be a multilinear function, $m \ge 0$, and suppose that for every bounded open set $U \subset Z$ and every $\epsilon > 0$ there is a compact set $C \subset U$ such that $\mass_{U \setminus C}(T) < \epsilon$. Then $\|T\|$ is a Borel regular outer measure that is concentrated on some $\sigma$-compact set.
\ep

\begin{proof}
It is clear that $\|T\|(\emptyset) = 0$ and that $\|T\|$ is monotone. To show that $\|T\|$ is $\sigma$-subadditive, let first $V_1,V_2,\ldots \subset Z$ be open, and put $V := \bigcup_{k=1}^\infty V_k$. Suppose $\bigl((f_\lambda,\pi_\lambda)\bigr)_{\lambda \in \Lambda}$ is a finite family as in the definition of $\mass_V(T)$. Choose a bounded open neighborhood $U \subset Z$ of $A := \bigcup_{\lambda \in \Lambda} \spt f_\lambda$, let $\epsilon > 0$, and let $C \subset U$ be a compact set with $\mass_{U \setminus C}(T) < \epsilon$. Put $K := C \cap A$ and $V_0 := Z \setminus K$. We have $K \subset \bigcup_{k=1}^\infty V_k$, thus by the compactness of $K$ there is an index $N$ such that $\bigcup_{k=1}^N V_k$ contains $K$. Furthermore, using the compactness of $K$ again, we see that the covering $(V_k)_{k=0,\dots,N}$ of $Z$ has a positive Lebesgue number. Then there exists a partition of unity $(\varrho_k)_{k=0,\dots,N}$ on $Z$ such that $\varrho_k \colon Z \to [0,1]$ is Lipschitz and $\spt \varrho_k \subset V_k$ for $k = 0,\dots,N$. For every $\lambda \in \Lambda$ we have $\spt(\varrho_0f_\lambda) \subset V_0 \cap A \subset U \setminus C$ and $\spt(\varrho_kf_\lambda) \subset V_k$ for $k = 1,\dots,N$; moreover $\sum_{\lambda \in \Lambda} |\varrho_k f_\lambda| \le 1$ for $k = 0,\dots,N$. We obtain
\[
 \sum_{\lambda \in \Lambda} T(f_\lambda,\pi_\lambda) 
 = \sum_{k=0}^N \sum_{\lambda \in \Lambda} T(\varrho_k f_\lambda,\pi_\lambda)
 \le \mass_{U \setminus C}(T) + \sum_{k=1}^N \mass_{V_k}(T) 
 < \epsilon + \sum_{k=1}^N \|T\|(V_k).
\]
It follows that $\|T\|(V) \le \sum_{k=1}^\infty \|T\|(V_k)$, and the same inequality for arbitrary sets $V_1,V_2,\dots$ is an immediate consequence. Thus $\|T\|$ is an outer measure. Whenever $A,B \subset Z$ with $d(A,B) := \inf\left\{d(x,y): x \in A, \,y \in B\right\} > 0$, then clearly $\|T\|(A \cup B) = \|T\|(A) + \|T\|(B)$. Hence, by Carath\'eodory's criterion, every Borel set is $\|T\|$-measurable, and by the definition of $\|T\|$ every set $A \subset Z$ is contained in a $G_\delta$ set $B$ with $\|T\|(B) = \|T\|(A)$. Thus $\|T\|$ is Borel regular. Writing $Z$ as the union of countably many bounded open sets $U_i$ and choosing a $\sigma$-compact set $\Sigma_i \subset U_i$ with $\|T\|(U_i \setminus \Sigma_i) = 0$ in each, we conclude that $\|T\|(Z \setminus \Sigma) = 0$ for $\Sigma := \bigcup_i\Sigma_i$, i.e.~$\|T\|$ is concentrated on the $\sigma$-compact set $\Sigma$.
\end{proof}

For a multilinear function $T$ satisfying the assumption of \propref{prop:sigma-cpt}, we define the support of $T$ as the closed set 
\begin{equation*}
 \spt T := \spt\|T\| 
 = \left\{z\in Z: \|T\|(\clB(z,r))>0\;\; \forall r>0\right\}.
\end{equation*}
If $\Sigma$ is a $\sigma$-compact set with $\|T\|(Z \setminus \Sigma) = 0$, 
then $\Sigma \setminus \spt T$ is contained in the union of countably many open balls with measure zero, thus
\begin{equation}\label{eq:spt}
 \|T\|(Z \setminus \spt T) = 0,
\end{equation}
i.e.~$\|T\|$ is concentrated on $\spt T$.

Now we return to metric functionals.

\bd \label{def:loc-fin-mass}
For $m \ge 0$, we denote by $\LWc{Z}{m}$ the vector space of all $m$-dimensional metric functionals $T$ on $Z$ (\defref{def:metric-fctl}) with the property that for every bounded open set $U \subset Z$ and every $\epsilon > 0$ there is a compact set $C \subset U$ such that $\mass_U(T) < \infty$ and $\mass_{U \setminus C}(T) < \epsilon$. Elements of $\LWc{Z}{m}$ will be called metric currents with locally finite mass.
\ed

By \propref{prop:sigma-cpt}, for every $T \in \LWc{Z}{m}$, $\|T\|$ is a Borel regular outer measure that is concentrated on some $\sigma$-compact set. The next result shows how $\|T\|$ controls $T$. We denote by $\Bor_Z$ the $\sigma$-algebra of Borel subsets of $Z$. By a Borel measure on $Z$ we mean a $\sigma$-additive function on $\Bor_Z$ with values in $[0,\infty]$.

\bp\label{prop:loc-fin-mass} 
Suppose that $T \in \LWc{Z}{m}$. Then
\begin{equation}\label{eq:int-ineq}
 |T(f,\pi_1,\dots,\pi_m)| 
 \le \prod_{i=1}^m\Lip(\pi_i\on{\spt f}) \int_Z|f|\, d\|T\|
\end{equation}
for all $(f,\pi_1,\dots,\pi_m) \in \form{Z}{m}$. Furthermore, if $\mu$ is a Borel measure on $Z$ such that $\mu$ is finite on bounded sets and this inequality holds with $\mu$ in place of $\|T\|$, then $\|T\| \le \mu$ on $\Bor_Z$.
\ep

\begin{proof}
We start with the case $m = 0$. Let $f \in \Lipbs(Z)$. There is no loss of generality in assuming $f \ge 0$. For any number $s \ge 0$, put $f_s := \min\{f,s\}$. Whenever $0 \le s < t$, we have $0 \le f_t - f_s \le t - s$ and hence $|T(f_t) - T(f_s)| = |T(f_t - f_s)| \le \|T\|(V)\,(t-s)$ for every open set $V$ containing the bounded set $\spt(f_t - f_s) = \overline{\{f > s\}}$; therefore
\[
 |T(f_t) - T(f_s)| \le \|T\|(\overline{\{f > s\}})\,(t-s).
\]
It follows that $s \mapsto T(f_s)$ is Lipschitz with constant $\le \|T\|(\spt f)$, hence almost everywhere differentiable, and $|(d/ds)T(f_s)| \le \|T\|(\overline{\{f > s\}})$ whenever the derivative exists. Since $T(f) = \int_0^\infty (d/ds) T(f_s) \,ds$, we conclude that
\[
 |T(f)| \le \int_0^\infty \biggl|\frac{d}{ds} T(f_s)\biggr| \,ds
 \le \int_0^\infty \|T\|(\overline{\{f > s\}}) \,ds = \int_Z f \,d\|T\|.
\]
This shows~\eqref{eq:int-ineq} in case $m = 0$. Now assume $m \ge 1$. Let first $(f,\pi) \in \Lipbs(Z) \times [\Lip_1(Z)]^m$, and consider the $0$-dimensional metric functional $T_\pi := T \rstr (1,\pi)$. Clearly $\|T_\pi\| \le \|T\|$, thus $T_\pi \in \LWc{Z}{0}$, and 
\[
 |T(f,\pi)| = |T_\pi(f)| \le \int_Z |f| \,d\|T_\pi\| \le \int_Z |f| \,d\|T\|.
\]
For a general $(f,\pi) \in \form{Z}{m}$, there exists $\bar\pi \in [\Lip(Z)]^m$ such that $\bar\pi = \pi$ on $\spt f$ and $\Lip(\bar\pi_i) = \Lip(\pi_i\on{\spt f})$ for $i = 1,\dots,m$. Then $T(f,\pi) = T(f,\bar\pi)$ by \lemref{lem:strict-loc}, and~\eqref{eq:int-ineq} follows.

As for the second assertion of the proposition, given such a measure $\mu$, we have $\mu(B) = \inf\left\{\mu(V): \text{$V \subset Z$ is open, $B \subset V$}\right\}$ for every Borel set $B \subset Z$ and $\mass_V(T) \le \mu(V)$ for every open set $V \subset Z$. This gives the result.
\end{proof}

Some basic examples of currents with locally finite mass are given as follows. Suppose $Z$ is a Lebesgue measurable subset of $\R^m$ with $\lm^m(\bdry Z) = 0$, and $\theta \colon Z \to \R$ is an $\lm^m$-measurable function such that $\int_{Z \cap U} |\theta| \,d\lm^m < \infty$ for every bounded open set $U \subset \R^m$. Then one obtains 
a current $[\theta] \in \LWc{Z}{m}$ by defining
\begin{equation}\label{eq:theta}
[\theta](f,\pi) := \int_Z \theta f \det(\nabla\pi)\,d\lm^m 
\end{equation}
for all $(f,\pi) \in \form{Z}{m}$ (cf.~\cite[Example~3.2]{AK-c} or~\cite[Proposition~2.6]{Lang-c}). It is not difficult to check that $\mass_V([\theta]) = 
\int_V |\theta|\,d\lm^m$ for every relatively open set $V \subset Z$. 

We conclude this section with some remarks regarding completeness of $Z$.
We did not impose a general completeness assumption on the underlying metric space, simply because this is not needed for the development of the theory. (The corresponding assumption in~\cite{AK-c} could equally well be avoided by some minor modifications.) However, the following simple example shows that some care is needed with incomplete spaces. Let $Z := (-\infty,0) \subset \R$ and $T := [1] \in \LWc{Z}{1}$, thus
\[
T(f,\pi) = \int_{-\infty}^0 f(s)\pi'(s) \,ds 
\]
for $(f,\pi) \in \Lipbs(Z) \times \Liploc(Z)$. As a ``constant'' current, $T$ should have no boundary in $Z$, however $\bdry T$ is the non-zero metric functional on $Z$ satisfying
\[
\bdry T(f) = \lim_{s \to 0-} f(s)
\]
for $f \in \Lipbs(Z)$. (In contrast, with the definitions from~\cite{Lang-c}, $\bdry T = 0$.) Furthermore $\mass(\bdry T) = 1$, yet $\bdry T \not\in \LWc{Z}{0}$ as $\|\bdry T\|(C) = 0$ for every compact set $C \subset Z$. Note that $\|\bdry T\|$ is not $\sigma$-subadditive in this case, and there is obviously no good definition of the support of $\bdry T$ in $Z$. Compare also~\eqref{eq:bdry-lip-mfd-current} in this regard, where $Z$ is a proper Lipschitz manifold.

\subsection{Extension to Borel functions}

We denote by $\Borbloc(Z)$ the algebra of all real-valued Borel functions on $Z$ that are bounded on bounded sets, and by $\Borb(Z)$ and $\Borbbs(Z)$ the subalgebras of bounded Borel functions and bounded Borel functions with bounded support, respectively. 

Due to~\eqref{eq:int-ineq}, every $T \in \LWc{Z}{m}$ naturally extends to a function 
\begin{equation*}
 T \colon \Borbbs(Z) \times [\Liploc(Z)]^m \to \R.
\end{equation*} 
In fact, whenever $f \in \Borbbs(Z)$ and $N$ is a bounded neighborhood of $\spt f$, there is a sequence $(f_k)$ in $\Lipbs(Z)$ such that $\spt f_k \subset N$ for all $k$ and $f_k \to f$ in $L^1(\|T\|)$. By~\eqref{eq:int-ineq}, $(T(f_k,\pi))$ is a Cauchy sequence for every $\pi \in [\Liploc(Z)]^m$, and $T(f,\pi)$ is declared as its limit, which is independent of the choice of $N$ and $(f_k)$. It is not difficult to show that the extended function $T$ satisfies the three conditions of \defref{def:metric-fctl} as well as inequality~\eqref{eq:int-ineq} with $\Borbbs(Z)$ in place of $\Lipbs(Z)$. The generalized inequality~\eqref{eq:int-ineq} also subsumes the strict locality property of \lemref{lem:strict-loc} for $f \in \Borbbs(Z)$. Furthermore, the extended functional is jointly continuous in all arguments in the following sense: If $(f,\pi),(f^j,\pi^j)\in\Borbbs(Z) \times [\Liploc(Z)]^m$, $(f^j,\pi^j) \to (f,\pi)$ pointwise everywhere as $j\to\infty$, $\sup_j\|f^j\|_\infty < \infty$, $\bigcup_j \spt f^j$ is bounded, and $\sup_{i,j}\Lip(\pi_i^j\on{A})<\infty$ for every bounded set $A\subset Z$, then
\begin{equation}\label{eq:joint-continuity}
 T(f^j,\pi^j) \to T(f,\pi)
\end{equation}
(cf.~\cite[Theorem~4.4]{Lang-c}). Due to~\eqref{eq:spt} and the generalized inequality~\eqref{eq:int-ineq}, the extended functional has the property that
\begin{equation}\label{eq:agree-on-spt}
 T(f,\pi) = T(f',\pi')
\end{equation}
whenever $f = f'$ and $\pi = \pi'$ on $\spt T$. From this it follows that $T$ may be viewed as an element of $\LWc{Y}{m}$ for any set $Y \subset Z$ containing $\spt T$ (cf.~\cite[Proposition~3.3]{Lang-c}). In particular, the push-forward $\varphi_\#T$ is still defined whenever $\varphi \colon D \to Z'$ is a map defined on a set $D \supset \spt T$ such that $\varphi\on{\spt T} \in \Liploc(\spt T,Z')$ and $\varphi^{-1}(A) \cap \spt T$ is bounded for every bounded set $A \subset Z'$.

The proof of the following product rule relies on the joint continuity property~\eqref{eq:joint-continuity}.

\bp\label{prop:prod-rule}
Let $T \in \LWc{Z}{m}$, $m \ge 1$. For all $f \in \Borbbs(Z)$ and $g,h,\pi_2,\dots,\pi_m\in \Liploc(Z)$,
\[
 T(f,gh,\pi_2,\dots,\pi_m) 
 = T(fg,h,\pi_2,\dots,\pi_m) + T(fh,g,\pi_2,\dots,\pi_m).
\]
\ep

\begin{proof}
This is shown as in~\cite[Proposition~2.4]{Lang-c}.  
\end{proof}

Let now $T \in \LWc{Z}{m}$ and $(g,\tau) \in \Borbloc(Z) \times [\Liploc(Z)]^l$ with $0 \le l \le m$. In view of the above extension, the restriction formula
\begin{equation*}
  T\rstr(g,\tau)\: (f,\pi) = T(fg,\tau,\pi)
\end{equation*}
for $(f,\pi) \in \form{Z}{m-l}$ remains meaningful and defines an $(m-l)$-dimensional metric functional $T\rstr(g,\tau)$. In fact, whenever $V \subset Z$ is an open set such that $\tau_j\on{V}$ is Lipschitz for $j = 1,\dots,l$, inequality~\eqref{eq:int-ineq} for the extended functional $T$ gives 
\begin{equation*}
 \mass_V(T\rstr(g,\tau)) \le \prod_{j=1}^l\Lip(\tau_j\on{V}) \int_V |g|\,d\|T\|,
\end{equation*}
which implies in particular that $T\rstr(g,\tau) \in \LWc{Z}{m-l}$. When $l = 1$ and $\tau = \tau_1$ is Lipschitz, this yields
\begin{equation}\label{eq:rstr-dtau}
 \|T\rstr(1,\tau)\|(A) \le \Lip(\tau) \,\|T\|(A)
\end{equation}
for every set $A \subset Z$. When $l = 0$, since $\|T\|(U) < \infty$ and $g\on{U}$ is bounded for every bounded set $U$, it follows that 
\begin{equation}\label{eq:rstr-g}
 \|T\rstr g\|(B) \le \int_B |g|\,d\|T\|
\end{equation}
for every Borel set $B \subset Z$. For the characteristic function $1_A$ of a Borel set $A \subset Z$, we write $T\rstr 1_A$ as $T\rstr A$. Then $\|T \rstr A\|(B) \le \|T\|(A \cap B) = (\|T\| \rstr A)(B)$ for every Borel set $B \subset Z$. In fact, since the same inequality holds for the complement $A^{\rm c}$, using the finiteness of $\|T\|$ on bounded sets and the identity $T = T\rstr A + T\rstr A^{\rm c}$ one easily concludes that
\begin{equation}\label{eq:rstr-a}
 \|T\rstr A\| = \|T\|\rstr A
\end{equation}
on $\Bor_Z$.

Let again $T \in \LWc{Z}{m}$, and let $\varphi\in\Liploc(Z,Z')$ be such that $\varphi^{-1}(A)$ is bounded whenever $A\subset Z'$ is. If $V \subset Z'$ is an open set such that $\varphi\on{\varphi^{-1}(V)}$ is $\lambda$-Lipschitz, \eqref{eq:int-ineq} yields
\begin{equation*}
 \mass_V(\varphi_{\#}T) 
 \le \lambda^m\,\|T\|(\varphi^{-1}(V))
 = \lambda^m\,(\varphi_\#\|T\|)(V). 
\end{equation*}
Given a bounded open set $U' \subset Z'$ and $\epsilon > 0$, there is a compact set $C \subset U := \varphi^{-1}(U')$ such that $\|T\|(U \setminus C) < \epsilon$, hence $C' := \varphi(C)$ is a compact subset of $U'$ with $(\varphi_\#\|T\|)(U' \setminus C')< \epsilon$. It follows that $\varphi_{\#}T \in \LWc{Z'}{m}$. Moreover, if $\varphi$ is Lipschitz, then $\mass_V(\varphi_{\#}T) \le \Lip(\varphi)^m\,(\varphi_\#\|T\|)(V)$ for every open set $V \subset Z'$, and since $\varphi_\#\|T\|$ is finite on bounded sets we have $(\varphi_\#\|T\|)(B) = \inf\left\{(\varphi_\#\|T\|)(V): \text{$V \subset Z'$ is open, $B \subset V$}\right\}$ for every Borel set $B \subset Z'$, so
\begin{equation*}
\|\varphi_\#T\| \le \Lip(\varphi)^m \,\varphi_\#\|T\|
\end{equation*}
on $\Bor_{Z'}$. 
We further note that the equation
\begin{equation}\label{eq:ext-push-fwd}
 \varphi_{\#}T\,(f,\pi) = T(f\circ\varphi, \pi\circ\varphi)
\end{equation}
remains valid for $(f,\pi) \in \Borbbs(Z') \times [\Liploc(Z')]^m$. To see this, choose a sequence $(f_k)$ in $\Lipbs(Z')$ such that $\bigcup_k\spt f_k$ is bounded and $f_k \to f$ in $L^1(\varphi_\#\|T\|)$, i.e.~$f_k \circ \varphi \to f \circ \varphi$ in $L^1(\|T\|)$. Then $f_k \to f$ in $L^1(\|\varphi_\#T\|)$ as well, and the result follows from inequality~\eqref{eq:int-ineq} for the extended functionals $\varphi_{\#}T$ and $T$.
Finally, we remark that if $(g,\tau) \in \Borbloc(Z') \times [\Liploc(Z')]^l$, $0 \le l \le m$, then $(g \circ \varphi,\tau \circ \varphi) \in \Borbloc(Z) \times [\Liploc(Z)]^l$ and 
\begin{equation*}
  \varphi_\#(T\rstr(g\circ\varphi,\tau\circ\varphi)) 
= (\varphi_\# T)\rstr(g,\tau),
 \end{equation*}
as is easily checked by means of~\eqref{eq:ext-push-fwd}. In particular,
\begin{equation}\label{eq:push-fwd-rstr}
\varphi_\#(T\rstr \varphi^{-1}(B)) = (\varphi_\#T)\rstr B
\end{equation}
for every Borel set $B \subset Z'$.

\subsection{Locally normal currents}

We now introduce the chain complex of locally normal currents.

\bd
For $m \ge 1$ we denote by $\LWnc{Z}{m}$ the vector space of all $T \in \LWc{Z}{m}$ with $\bdry T\in\LWc{Z}{m-1}$, and we put $\LWnc{Z}{0} := \LWc{Z}{0}$. Elements of $\LWnc{Z}{m}$ will be called locally normal currents.
\ed

Let $m \ge 1$ and $g \in \Liploc(Z)$, and suppose first that $T \in \LWc{Z}{m}$. For $(f,\pi) \in \form{Z}{m-1}$ and $\sigma \in \Lipbs(Z)$ with $\sigma\on{\spt f} = 1$, \propref{prop:prod-rule} gives $T(\sigma,fg,\pi) = T(\sigma g,f,\pi) + T(f,g,\pi)$, which corresponds to the identity
\begin{equation}\label{eq:bdry-rstr}
(\bdry T)\rstr g = \bdry(T\rstr g) + T \rstr (1,g)
\end{equation}
of $(m-1)$-dimensional metric functionals. Since $T \rstr (1,g) \in \LWc{Z}{m-1}$, it follows that $(\bdry T)\rstr g\in\LWc{Z}{m-1}$ if and only if $\bdry(T\rstr g)\in\LWc{Z}{m-1}$. Now let $T \in \LWnc{Z}{m}$. Then $(\bdry T)\rstr g\in\LWc{Z}{m-1}$ and hence $T \rstr g \in \LWnc{Z}{m}$. Furthermore, if $g$ is Lipschitz, combining~\eqref{eq:bdry-rstr} with~\eqref{eq:rstr-dtau} and~\eqref{eq:rstr-g} we see that
\begin{equation*}
\|\bdry(T\rstr g)\|(B) \le \Lip(g) \,\|T\|(B) + \int_B|g|\,d\|\bdry T\| 
\end{equation*}
for every Borel set $B \subset Z$. Push-forwards of locally normal currents are locally normal. The following simple criterion will be useful:

\bl\label{lem:char-LWnc}
Suppose $T \colon \form{Z}{m} \to \R$ is a function, $(\sigma_i)$ is a sequence in $\Lipbs(Z)$ such that every bounded set $A \subset Z$ is contained in $\{\sigma_i = 1\}$ for some $i$, and $T \rstr \sigma_i \in \LWnc{Z}{m}$ for every $i$. Then $T \in \LWnc{Z}{m}$.
\el

\begin{proof}
It is easily checked that $T$ is a metric functional and that $\mass_V(T) = \mass_V(T \rstr \sigma_i)$ whenever $V$ is a bounded open set and $\sigma_i\on{V} = 1$. Thus $T \in \LWc{Z}{m}$. Moreover, in case $m \ge 1$, $\bdry T(f,\pi) = \bdry(T \rstr \sigma_i)(f,\pi)$ whenever $(f,\pi) \in \form{Z}{m-1}$ and $\sigma_i\on{\spt f} = 1$, hence also $\mass_V(\bdry T) = \mass_V(\bdry(T \rstr \sigma_i))$.
\end{proof}

\subsection{Relation to Ambrosio--Kirchheim currents}

We now discuss the relation to the theory of Ambrosio--Kirchheim. We briefly recall that a current $T \in \AKc_m(Z)$ in the sense of \cite{AK-c} is a multilinear function $T \colon \Lipb(Z) \times [\Lip(Z)]^m \to \R$ such that 
\begin{equation}\label{eq:AK-cont}
T(f,\pi^j) \to T(f,\pi)
\end{equation}
whenever $\pi^j \to \pi$ pointwise with $\sup_{i,j}\Lip(\pi^j_i) < \infty$;
furthermore, by assumption, there exists a finite Borel measure $\mu$ on $Z$ such that 
\begin{equation}\label{eq:AK-int-ineq}
|T(f,\pi)| \le \prod_{i=1}^m \Lip(\pi_i\on{\spt f}) \int_Z |f| \,d\mu
\end{equation}
for all $(f,\pi) \in \Lipb(Z) \times [\Lip(Z)]^m$ (in particular $T(f,\pi) = 0$ if some $\pi_i$ is constant on $\spt f$). There is a least Borel measure with this property, denoted $\|T\|$, and there exists a $\sigma$-compact set $\Sigma \subset Z$ such that $\|T\|(Z \setminus \Sigma) = 0$ (cf.~Lemma~2.9 in \cite{AK-c} and the remark thereafter; note also that the proof of this lemma requires completeness of the underlying metric space). As above, $\mass(T) := \|T\|(Z)$, $\spt T := \spt\|T\|$, and $\|T\|(Z \setminus \spt T) = 0$. For $m \ge 1$, the functional $\bdry T$ is defined by $\bdry T(f,\pi_1,\dots,\pi_{m-1}) = T(1,f,\pi_1,\dots,\pi_{m-1})$, $\AKnc_m(Z) := \left\{T \in \AKc_m(Z): \bdry T \in \AKc_{m-1}(Z)\right\}$, and $\AKnc_0(Z) := \AKc_0(Z)$. Every $T \in \AKc_m(Z)$ extends to a multilinear function $T \colon \Borb(Z) \times [\Lip(Z)]^m \to \R$ such that~\eqref{eq:AK-cont} and~\eqref{eq:AK-int-ineq} still hold for $f \in \Borb(Z)$. Given $T \in \AKc_m(Z)$, $(g,\tau) \in \Borb(Z) \times [\Lip(Z)]^l$ with $0 \le l \le m$, and $\varphi \in \Lip(Z,Z')$, the restriction $T \rstr (g,\tau) \in \AKc_{m-l}(Z)$ and the push-forward $\varphi_\#T \in \AKc_m(Z')$ are defined in the same way as in our case. 

Let now $T \in \LWc{Z}{m}$ and $g \in \Borbbs(Z)$. Since the support of $g$ is bounded, the formula
\[
T\rstr g\: (f,\pi) = T(f g,\pi)
\] 
remains meaningful for $(f,\pi) \in \Lipb(Z) \times [\Lip(Z)]^m$. Thus, the restriction $T\rstr g \in \LWc{Z}{m}$ may as well be viewed as a function on $\Lipb(Z) \times [\Lip(Z)]^m$, which we denote $T_g$ for the moment. It is easily checked that $T_g$ is an element of $\AKc_m(Z)$: Clearly $T_g$ is multilinear and satisfies~\eqref{eq:AK-cont}; furthermore, by~\eqref{eq:rstr-g}, $\|T \rstr g\|$ is concentrated on the bounded set $\spt g$, it thus follows from \propref{prop:loc-fin-mass} for $T \rstr g$ that~\eqref{eq:AK-int-ineq} holds for $T_g$ and that $\|T_g\| = \|T \rstr g\|$ on $\Bor_Z$. In case $m \ge 1$ and $g \in \Lipbs(Z)$, when $T \rstr g \in \LWnc{Z}{m}$, we have $T_g \in \AKnc_m(Z)$ and $\|\bdry(T_g)\| = \|\bdry(T \rstr g)\|$. To see this, choose $\sigma \in \Lipbs(Z)$ with $\sigma\on{\spt g} = 1$; then $(\bdry(T \rstr g))_\sigma \in \AKc_{m-1}(Z)$ and $\|(\bdry(T \rstr g))_\sigma\| = \|\bdry(T \rstr g) \rstr \sigma\|$, and it is not difficult to verify that $(\bdry(T \rstr g))_\sigma = \bdry(T_g)$ and $\bdry(T \rstr g) \rstr \sigma = \bdry(T \rstr g)$. From now on we write again $T \rstr g$ instead of $T_g$; an expression like $T \rstr g \in \AKc_m(Z)$ will indicate that a function on $\Lipb(Z) \times [\Lip(Z)]^m$ is understood.

We show next that every $T\in\LWc{Z}{m}$ with finite mass 
determines an Ambrosio--Kirchheim current $T' \in \AKc_m(Z)$.

\bp
Let $T\in\LWc{Z}{m}$ with $\mass(T) < \infty$. Then there exists $T'\in\AKc_m(Z)$ such that, whenever $(\sigma_n)$ is a sequence in $\Lipbs(Z)$ such that $|\sigma_n|\le 1$ and $\sigma_n \to 1$ uniformly on bounded sets, the restrictions $T\rstr\sigma_n \in \AKc_m(Z)$ converge in mass to $T'$.
\ep

\begin{proof}
Let $A\subset Z$ be a bounded Borel set, and let $0 < \varepsilon < 1$. Suppose $\varrho\in\Lipbs(Z)$, $|\varrho|\leq 2$, and $|\varrho|\leq\varepsilon$ on $A$. Using~\eqref{eq:rstr-g} we obtain
\begin{equation*}
 \mass(T\rstr\varrho) = \|T \rstr\varrho\|(Z) \le \int_Z |\varrho|\,d\|T\| \le \varepsilon \|T\|(A) + 2\|T\|(A^{\rm c}).
\end{equation*}
In particular, if $\|T\|(A^{\rm c})\leq\varepsilon\|T\|(Z)$ and if $\sigma, \sigma'\in\Lipbs(Z)$ are such that $|\sigma|, |\sigma'|\leq 1$ and $\sigma\on{A}, \sigma'\on{A}\geq 1-\varepsilon$, then
\begin{equation*}
 \mass(T\rstr\sigma - T\rstr\sigma') = \mass(T \rstr (\sigma- \sigma'))
\leq 3\varepsilon\|T\|(Z).
\end{equation*}
Since $\AKc_m(Z)$ is $\mass$-complete, the result follows.
\end{proof}

Conversely, given an Ambrosio--Kirchheim current $T' \in \AKc_m(Z)$, one obtains a well-defined element $T \in \LWc{Z}{m}$ by putting 
\begin{equation*}
T(f,\pi) := T'(f,\pi')
\end{equation*}
for $(f,\pi) \in \form{Z}{m}$ and any $\pi' \in [\Lip(Z)]^m$ with $\pi_i'\on{\spt f} = \pi_i\on{\spt f}$; moreover $\|T\| = \|T'\|$ on $\Bor_Z$. In case $m \ge 1$, it follows that $\bdry T(f,\pi) = \bdry T'(f,\pi')$ whenever $(f,\pi) \in \form{Z}{m-1}$ and $\pi' \in [\Lip(Z)]^{m-1}$ with $\pi_i'\on{\spt f} = \pi_i\on{\spt f}$. In particular, if $T' \in \AKnc_m(Z)$, then $T \in \LWnc{Z}{m}$ and $\|\bdry T\| = \|\bdry T'\|$ on $\Bor_Z$. 

\subsection{Slices}

Let $T\in\LWnc{Z}{m}$, $m\geq1$, and let $\varrho \in \Liploc(Z)$. For every $r \in \R$ we define the $(m-1)$-dimensional metric functional
\begin{equation*}
\langle T, \varrho, r\rangle 
:= \bdry(T\rstr\{\varrho\leq r\}) - (\bdry T)\rstr\{\varrho\leq r\},
 \end{equation*}
called the (right-handed) slice of $T$ at $r$, with respect to $\varrho$. For every $g \in \Liploc(Z)$, the slice of $T \rstr g \in \LWnc{Z}{m}$ at $r$ is given by
\begin{equation}\label{eq:slice-rstr}
\langle T \rstr g, \varrho, r\rangle = \langle T, \varrho, r\rangle \rstr g.
\end{equation}
To see this, put $T_r := T\rstr\{\varrho\leq r\} \in \LWc{Z}{m}$ and $(\bdry T)_r := (\bdry T) \rstr \{\varrho\leq r\} \in \LWc{Z}{m-1}$, so that $\langle T, \varrho, r\rangle \rstr g = (\bdry T_r) \rstr g - (\bdry T)_r \rstr g$. By~\eqref{eq:bdry-rstr},
\[
\bdry((T \rstr g) \rstr \{\varrho\leq r\}) = \bdry(T_r \rstr g) 
= (\bdry T_r) \rstr g - T_r \rstr (1,g)
\]
and $\bdry(T \rstr g) = (\bdry T) \rstr g - T \rstr (1,g)$, hence 
\[
(\bdry(T \rstr g)) \rstr \{\varrho\leq r\} = (\bdry T)_r \rstr g - T_r \rstr (1,g),
\]
and~\eqref{eq:slice-rstr} follows. As in the theories in \cite{AK-c} and \cite{Lang-c} we have:

\bt\label{thm:slice-mass-bound}
Let $T\in\LWnc{Z}{m}$, $m \ge 1$, and let $\varrho\in\Lip(Z)$. Then for almost every $r\in\R$ the slice $\langle T, \varrho, r\rangle$ is a locally normal current with 
 \begin{equation*}
   \spt \langle T, \varrho, r\rangle \subset \{\varrho = r\} \cap \spt T.
 \end{equation*}
 Moreover, for every Borel set $A\subset Z$ and for $-\infty<r_0< r_1<\infty$,
 \begin{equation*}
  \int_{r_0}^{r_1}\|\langle T, \varrho, r\rangle \|(A)\,dr \leq \Lip(\varrho)\,\|T\|(A\cap\{r_0 < \varrho < r_1\}).
 \end{equation*} 
\et

\begin{proof}
We choose a sequence $(\sigma_i)$ in $\Lipbs(Z)$ such that every bounded set $A \subset Z$ is contained in $\{\sigma_i=1\}$ for some $i$. For every $i$, $T_i := T \rstr \sigma_i \in \AKnc_m(Z)$, and $\spt T_i \subset \spt T$. Applying the slicing theorem~\cite[Theorem 5.6]{AK-c} to each $T_i$, we conclude that there is a set $N \subset \R$ of measure zero such that for every $r \in \R \setminus N$, 
\[
\langle T_i,\varrho,r\rangle = \bdry(T_i\rstr\{\varrho\leq r\}) - (\bdry T_i)\rstr\{\varrho\leq r\} \in \AKnc_{m-1}(Z)
\]
for all $i$, $\spt(\langle T_i,\varrho,r\rangle) \subset \{\varrho = r\} \cap \spt T_i$, and 
\begin{equation*}
  \int_{r_0}^{r_1}\|\langle T_i, \varrho, r\rangle \|(A)\,dr \leq \Lip(\varrho)\,\|T_i\|(A\cap\{r_0 < \varrho < r_1\})
\end{equation*} 
for every Borel set $A\subset Z$ and for $-\infty<r_0< r_1<\infty$. Now we view $\langle T_i,\varrho,r\rangle$ as an element of $\LWnc{Z}{m-1}$. It follows from~\eqref{eq:slice-rstr} that $\langle T,\varrho,r\rangle \rstr \sigma_i$ and $\langle T_i,\varrho,r\rangle$ agree as functions on $\form{Z}{m-1}$. Thus, for every $r \in \R \setminus N$, $\langle T,\varrho,r\rangle$ is a metric functional with the property that $\langle T,\varrho,r\rangle \rstr \sigma_i \in \LWnc{Z}{m-1}$ for all $i$. Therefore $\langle T,\varrho,r\rangle \in \LWnc{Z}{m-1}$ by \lemref{lem:char-LWnc}. If $A \subset Z$ is a bounded Borel set and $i$ is such that $\sigma_i\on{A} = 1$, then 
\begin{equation*}
\|\langle T,\varrho,r\rangle\|(A) = \|\langle T,\varrho,r\rangle \rstr \sigma_i\|(A) = \|\langle T_i, \varrho, r\rangle \|(A)
\end{equation*}
and $\|T_i\|(A) = \|T\|(A)$. We conclude that $\|\langle T,\varrho,r\rangle\|$ is concentrated on $\{\varrho = r\} \cap \spt T$ and that the claimed inequality holds for bounded Borel sets, hence also for arbitrary Borel sets $A \subset Z$.
\end{proof}

\subsection{Locally integer rectifiable and integral currents}

We call a subset of $Z$ a compact $m$-rectifiable set if it is the union of finitely many sets of the form $\lambda(K)$ where $K \subset \R^m$ is compact and $\lambda \in \Lip(K,Z)$. A compact $0$-rectifiable set is just a finite set. For condition~(ii) below we recall the basic examples of currents defined in~\eqref{eq:theta}.

\bd\label{def:LWirc}
For $m \ge 0$, we denote by $\LWirc{Z}{m}$ the set of all $m$-dimensional metric functionals on $Z$ with the following two properties:
\begin{enumerate}
\item
For every bounded open set $U \subset Z$ and every $\epsilon > 0$ there is a compact $m$-rectifiable set $C \subset U$ such that $\mass_U(Z) < \infty$ and $\mass_{U \setminus C}(T) < \epsilon$, in particular $T \in \LWc{Z}{m}$;
\item
for every bounded Borel set $B \subset Z$ and every $\pi \in \Lip(Z,\R^m)$ there exists $\theta \in L^1(\R^m,\Z)$ such that $\pi_\#(T \rstr B) = [\theta]$.
\end{enumerate}
Elements of $\LWirc{Z}{m}$ are called locally integer rectifiable currents.
\ed

By~(i), $T \in \LWc{Z}{m}$, and $\|T\|$ is concentrated on the union of countably many sets of the form $\lambda(K)$ as above. Conversely, if $T \in \LWc{Z}{m}$ and $\|T\|$ is concentrated on such a union, then clearly $T$ satisfies~(i). In~(ii), $\pi_\#(T \rstr B)$ is defined as an element of $\LWc{\R^m}{m}$ according to the remark after~\eqref{eq:agree-on-spt}, in fact 
\[
\pi_\#(T \rstr B)(f,g) = T(1_B(f \circ \pi),g \circ \pi)
\] 
for every $(f,g) \in \Borbbs(\R^m) \times [\Liploc(\R^m)]^m$. We also remark that it suffices to verify condition~(ii) for bounded open sets $B \subset Z$, cf.~the proof of~\cite[Lemma~7.3]{Lang-c}. In case $m = 0$, an element $T \in \LWirc{Z}{0}$ is just a function $T \colon \Lipbs(Z) \to \R$ of the following form: There exist a set $E \subset Z$ and a function $\theta \colon E \to \Z$ such that every bounded subset of $E$ is finite and  
\[
T(f) = \sum_{z \in E} \theta(z) f(z)
\]
for every $f \in \Lipbs(Z)$. Clearly $\LWirc{Z}{m}$ forms an additive abelian group. Let $T \in \LWirc{Z}{m}$. If $\varphi\in\Liploc(Z,Z')$ is such that $\varphi^{-1}(A)$ is bounded whenever $A \subset Z'$ is, then $\varphi_\#T \in \LWirc{Z'}{m}$; this uses~\eqref{eq:push-fwd-rstr} and the fact that $\pi_\# \circ \varphi_\# = (\pi \circ \varphi)_\#$. If $A \subset Z$ is a Borel set, then obviously $T \rstr A \in \LWirc{Z}{m}$.

\bp
If a current $T \in \LWc{Z}{m}$ satisfies condition~{\rm (ii)} of \defref{def:LWirc} with $L^1(\R^m)$ in place of $L^1(\R^m,\Z)$, then $\|T\|$ is absolutely continuous with respect to $\hm^m$.
\ep

\begin{proof}
It suffices to show that $\|T\|(C) = 0$ for every bounded closed set $C \subset Z$ with $\hm^m(C) = 0$. Suppose to the contrary that there is such a set $C$ with $\|T\|(C) > 0$. By~\eqref{eq:rstr-a}, $\mass(T \rstr C) > 0$, hence there exists $(f,\pi) \in \Lipbs(Z) \times [\Lip(Z)]^m$ such that $T(1_Cf,\pi) = T \rstr C\:(f,\pi) \ne 0$. Approximating $f$ by simple functions, and using the continuity of the extended functional in the first argument, we find a closed set $B \subset C$ such that $T(1_B,\pi) \ne 0$. Since $\lm^m(\pi(B)) = 0$, there is a bounded Borel set $N \subset \R^m$ such that $\pi(B) \subset N$ and $\lm^m(N) = 0$. Now
\[
T(1_B,\pi) = T(1_B(1_N \circ \pi),\pi) 
= \pi_\#(T \rstr B)(1_N,\id) = [\theta](1_N,\id)
\]
for some $\theta \in L^1(\R^m)$. Since $\lm^m(N) = 0$, $[\theta](1_N,\id) = 0$, a contradiction. 
\end{proof}

We now introduce the chain complex of locally integral currents.

\bd
For $m \ge 1$ we denote by $\LWic{Z}{m}$ the abelian group of all $T \in \LWirc{Z}{m}$ with $\bdry T\in\LWirc{Z}{m-1}$, and we put $\LWic{Z}{0} := \LWirc{Z}{0}$. Elements of $\LWic{Z}{m}$ will be called locally integral currents.
\ed

In particular, locally integral currents are locally normal. In fact, \thmref{thm:bdry-rect} below will show that $\LWic{Z}{m} = \LWirc{Z}{m} \cap \LWnc{Z}{m}$. In analogy with \lemref{lem:char-LWnc} we have:

\bl\label{lem:char-LWic}
Suppose $T \in \LWc{Z}{m}$, $(A_i)$ is a sequence of bounded Borel subsets of $Z$ such that every bounded set $A \subset Z$ is contained in some $A_i$, and $T \rstr A_i \in \LWic{Z}{m}$ for every $i$. Then $T \in \LWic{Z}{m}$.
\el

\begin{proof}
It is easily checked that $T \in \LWirc{Z}{m}$. Moreover, in case $m \ge 1$, it follows from the strict locality of the extended functional $T$ that $\bdry T(f,\pi) = \bdry(T \rstr A_i)(f,\pi)$ whenever $(f,\pi) \in \form{Z}{m-1}$ and $\{f \ne 0\} \subset A_i$, and this yields $\bdry T \in \LWirc{Z}{m-1}$.
\end{proof}

Next, we deduce the Boundary Rectifiability Theorem for locally integer rectifiable currents from the corresponding result in~\cite{AK-c}. We denote by $\AKirc_m(Z)$ and $\AKic_m(Z)$ the spaces of integer rectifiable and integral currents in $Z$, as defined in~\cite[Definition~4.2]{AK-c}. 

\bt\label{thm:bdry-rect}
If $T\in\LWirc{Z}{m}$, $m\ge 1$, and $\bdry T \in \LWc{Z}{m-1}$, then $\bdry T\in\LWirc{Z}{m-1}$, i.e.~$T \in \LWic{Z}{m}$.
\et

\begin{proof}
Note that $T \in \LWnc{Z}{m}$. Let $\varrho$ be the distance function to a fixed point $z_0 \in Z$, choose a sequence $0 < r_1 < r_2 < \ldots \to \infty$ such that $\langle T,\varrho,r_i\rangle \in \LWnc{Z}{m-1}$ for all $i$, and put $A_i := \clB(z_0,r_i)$. Then $\bdry(T\rstr A_i) = \langle T,\varrho,r_i\rangle + (\bdry T) \rstr A_i \in \LWc{Z}{m-1}$, thus $T \rstr A_i \in \LWirc{Z}{m} \cap \LWnc{Z}{m}$. Now we view $T \rstr A_i$ as an element of $\AKnc_m(Z)$. Then clearly $T \rstr A_i$ also belongs to $\AKirc_m(Z)$. By~\cite[Theorem~8.6]{AK-c}, $\bdry(T \rstr A_i) \in \AKirc_{m-1}(Z)$. Interpreting $\bdry(T \rstr A_i)$ again as element of $\LWc{Z}{m-1}$, we conclude that $\bdry(T \rstr A_i) \in \LWirc{Z}{m-1}$. As this holds for every $A_i$, we have $T \in \LWic{Z}{m}$ by \lemref{lem:char-LWic}.
\end{proof}

As a consequence, one obtains the following supplement to \thmref{thm:slice-mass-bound}: Whenever $T\in\LWic{Z}{m}$, $m \ge 1$, $\varrho\in\Liploc(Z)$, and $\langle T, \varrho, r\rangle\in\LWnc{Z}{m-1}$ for some $r \in \R$, then $\bdry(T\rstr \{\varrho \le r\}) \in \LWc{Z}{m-1}$, hence 
\begin{equation*}
T \rstr \{\varrho \le r\} \in \LWic{Z}{m}
\end{equation*} 
and $\langle T, \varrho, r\rangle\in\LWic{Z}{m-1}$.

Finally, we deduce the Closure Theorem for locally integral currents from the corresponding result in~\cite{AK-c}.

\bt
Suppose $m \ge 1$, and $(T_n)$ is a sequence in $\LWic{Z}{m}$ that converges weakly to some $T \in \LWnc{Z}{m}$, with
 \begin{equation*}
  \sup_n[\|T_n\|(A) + \|\bdry T_n\|(A)] < \infty
 \end{equation*}
 for every bounded Borel set $A\subset Z$. Then $T\in\LWic{Z}{m}$.
\et

\begin{proof}
Let $\varrho$ be the distance function to a fixed point $z_0 \in Z$. As in the proofs of~\cite[Proposition~8.3]{AK-c} and~\cite[Proposition~6.6]{Lang-c} one shows that for almost every $r > 0$ there exists a subsequence $(n(k))$ such that $\langle T_{n(k)},\varrho,r\rangle \in \LWnc{Z}{m-1}$ for all $k$,
\[
\sup_k\mass(\langle T_{n(k)},\varrho,r\rangle) < \infty,
\]
and $T_{n(k)} \rstr A \to T \rstr A$ weakly, where $A := \clB(z_0,r)$. It follows that $T_{n(k)} \rstr A \in \LWic{Z}{m}$ and 
\[
\sup_k[\mass(T_{n(k)} \rstr A) + \mass(\bdry(T_{n(k)} \rstr A))] < \infty.
\]
In addition, for almost every $r \in \R$, $\langle T,\varrho,r\rangle \in \LWnc{Z}{m-1}$ and hence $T \rstr A \in \LWnc{Z}{m}$. Now we interpret $T_{n(k)} \rstr A$ and $T \rstr A$ as elements of $\AKic_m(Z)$ and $\AKnc_m(Z)$, respectively. By~\cite[Theorem~8.5]{AK-c}, $T \rstr A \in \AKic_{m}(Z)$, hence $T \rstr A \in \LWic{Z}{m}$ as a function on $\form{Z}{m}$. In view of \lemref{lem:char-LWic}, the result follows.
\end{proof}

\subsection{Manifolds as currents}\label{subsection:manifolds}

Every connected and oriented Riemannian manifold $M$ that is complete as a metric space gives rise to a locally integral current $[M]$ in $M$ of the same dimension. The same is true for proper, oriented Lip\-schitz manifolds, as we show now. Recall that a metric space $Z$ is an $m$-dimensional Lipschitz manifold if it can be covered by charts $(U_\alpha,\varphi_\alpha)$ where $U_\alpha \subset Z$ is open and $\varphi_\alpha$ is a bi-Lipschitz map from $U_\alpha$ onto a relatively open subset of $H_\alpha := \{\lambda_\alpha \ge 0\}$, for some linear function $\lambda_\alpha\colon \R^m \to \R$. If $m \ge 1$, the boundary $\bdry Z$ is the $(m-1)$-dimensional Lipschitz manifold consisting of all $z \in Z$ such that $\varphi_\alpha(z) \in \bdry H_\alpha$ for some $\alpha$. A Lipschitz manifold $Z$ of dimension $m \ge 1$ is said to be orientable if it admits an atlas $\{(U_\alpha, \varphi_\alpha)\}_{\alpha\in A}$ such that $\det(\nabla(\varphi_\alpha\circ\varphi_\beta^{-1}))>0$ almost everywhere on $\varphi_\beta(U_\alpha\cap U_\beta)$, for all $\alpha, \beta\in A$. An orientation is a maximal such atlas. If $m\geq 2$, then an orientation on $Z$ induces an orientation on~$\bdry Z$.

Let now $Z$ be a proper, oriented, $m$-dimensional Lipschitz manifold. Choose a locally finite (hence countable) atlas $\{(U_\alpha, \varphi_\alpha)\}_{\alpha\in A}$ of positively oriented charts. Let furthermore $(\varrho_\alpha)$ be a locally Lipschitz partition of unity on $Z$ with $\spt \varrho_\alpha \subset U_\alpha$. For $(f,\pi)\in\form{Z}{m}$ we define
\begin{equation*}
 \begin{split}
 [Z](f, \pi)&:= \sum_{\alpha \in A}(\varphi^{-1}_\alpha)_{\#}[\varrho_\alpha\circ\varphi^{-1}_\alpha](f, \pi)\\
  &\phantom{:}= \sum_{\alpha \in A} \int_{\varphi_\alpha(U_\alpha)}((\varrho_\alpha f)\circ\varphi_\alpha^{-1})\det\left(\nabla(\pi\circ\varphi_\alpha^{-1})\right) \,d\lm^m.
 \end{split}
\end{equation*}
Since $\spt f$ is compact and the chosen atlas locally finite, only finitely many terms in these sums are non-zero. Furthermore we clearly have $[Z]\in\LWc{Z}{m}$ because $(\varphi^{-1}_\alpha)_{\#}[\varrho_\alpha\circ\varphi^{-1}_\alpha]\in\LWc{Z}{m}$ for every $\alpha$ and because the atlas is locally finite. It follows from the lemma below that $[Z]\in\LWirc{Z}{m}$ and that $[Z]$ is independent of the particular choices of atlas and partition of unity.

\bl \label{lem:z-rstr-g}
Let $(U,\psi)$ be a positively oriented chart of $Z$, and let $g \in \Borbloc(Z)$ with $\spt g \subset U$. Then
\begin{equation*}
[Z]\rstr g = \psi^{-1}_\#[g \circ \psi^{-1}].
\end{equation*}
\el

\begin{proof}
For $(f,\pi) \in \form{\psi(U)}{m}$, we have
\begin{equation*}
\begin{split}
\psi_\#([Z]\rstr g)(f,\pi) 
&= \sum_{\alpha \in A} \int_{\varphi_\alpha(U_\alpha \cap U)}((\varrho_\alpha g)\circ \varphi_\alpha^{-1})(f \circ\psi \circ \varphi_\alpha^{-1})\det\left(\nabla(\pi\circ\psi\circ\varphi_\alpha^{-1})\right)\,d\lm^m\\
&= \sum_{\alpha \in A} \int_{\psi(U_\alpha \cap U)}((\varrho_\alpha g)\circ\psi^{-1})f\det(\nabla \pi) \,d\lm^m\\
&= [g \circ \psi^{-1}](f,\pi).
\end{split}
\end{equation*}
This proves the lemma.
\end{proof}

We now show that if $m\geq 2$, then 
\begin{equation}\label{eq:bdry-lip-mfd-current}
 \bdry[Z] = [\bdry Z].
\end{equation} 
Let $(f,\pi)\in\form{Z}{m-1}$, and choose $\sigma_\alpha \in \Lipbs(Z)$ such that $\spt \sigma_\alpha \subset U_\alpha$ and $\sigma_\alpha = 1$ on $\spt(\varrho_\alpha f)$. Then
\[
\bdry[Z](f,\pi) = \sum_{\alpha \in A}\bdry[Z](\varrho_\alpha f,\pi)
= \sum_{\alpha \in A}[Z](\sigma_\alpha,\varrho_\alpha f,\pi).
\]
Furthermore, by \lemref{lem:z-rstr-g} and Stokes' theorem, generalized to Lipschitz functions by bounded smooth approximation,
\begin{align*}
[Z](\sigma_\alpha,\varrho_\alpha f,\pi) 
  &= 
\int_{\varphi_\alpha(U_\alpha)}\det\left(\nabla((\varrho_\alpha f)\circ\varphi_\alpha^{-1}, \pi\circ\varphi_\alpha^{-1})\right)\,d\lm^m\\
  &= 
\int_{\varphi_\alpha(U_\alpha)} d((\varrho_\alpha f)\circ\varphi_\alpha^{-1})\wedge d(\pi_1\circ\varphi_\alpha^{-1})\wedge\ldots\wedge d(\pi_{m-1}\circ\varphi_\alpha^{-1})\\
  &= 
\int_{\bdry H_\alpha \cap \varphi_\alpha(U_\alpha)}((\varrho_\alpha f)\circ\varphi_\alpha^{-1})\,d(\pi_1\circ\varphi_\alpha^{-1})\wedge\ldots\wedge d(\pi_{m-1}\circ\varphi_\alpha^{-1}).
\end{align*}
This gives~\eqref{eq:bdry-lip-mfd-current}. In particular, if $m \ge 2$,
$[Z]$ is a locally integral current, and it is not difficult to check that 
this is true also when $m=1$. 


\section{Proofs of the main results}\label{Section:main-results}

\subsection{The pointed compactness theorem}\label{subsection:ptd-cptness}

We now turn to our main result, \thmref{thm:main}, whose proof relies on the arguments of \cite{Wenger-cpt}. The proposition below summarizes some key facts established in Lemma~5.1 and the first part of the proof of Theorem~1.2 in that paper. For $n \in \N$, let $X_n$ be a complete metric space with basepoint $x_n$, and let $m \in \N$. Replacing $X_n$ by $l^\infty(X_n)$ if necessary, we may assume by~\cite{Wenger-isopineq} that for $k = 1,\ldots,m$, $X_n$ admits an isoperimetric inequality of Euclidean type for $\AKic_k(X_n)$ with constant $D_k$. This means that for every $R \in \AKic_k(X_n)$ with $\bdry R = 0$ there exists $S \in \AKic_{k+1}(X_n)$ with $\bdry S = R$ such that 
\[
\mass(S) \le D_k \mass(R)^{(k+1)/k}.
\]
Furthermore, since every closed ball in $l^\infty(X_n)$ is a $1$-Lipschitz retract, we may assume that $\spt S \subset B$ whenever $\spt R$ is contained
in some fixed closed ball $B \subset X_n$.
Now, fix integers $1=j_1<j_2<j_3<\ldots$ and positive numbers $\frac{1}{2}>\delta_1>\delta_2>\ldots$ such that
\begin{equation*}
 \Delta:=\sum_{i=1}^\infty \delta_i<\infty.
\end{equation*}
Recall that a sequence of compact metric spaces $K_n$ is said to be uniformly compact if the diameters are uniformly bounded and if for every $\epsilon > 0$ there is $N(\epsilon) \in \N$ such that every $K_n$ can be covered by $N(\epsilon)$ open balls of radius $\epsilon$.

\bp\label{prop:good-decomp-summary}
 Let $R, C>0$ and suppose that for every $n \in \N$, $T_n\in\AKic_m(X_n)$ satisfies $\spt T_n\subset \clB(x_n, R)$ and
 \begin{equation*}
  \mass(T_n) + \mass(\bdry T_n) \leq C.
 \end{equation*}
 Then there exist currents $T_n^1, \dots, T_n^{j_n+1},U_n^1, \dots, U_n^{j_n+1} \in \AKic_m(X_n)$ with support in $\clB(x_n, R)$ such that 
 $$T_n = T_n^1 + \ldots + T_n^{j_n+1} + U_n^1 + \ldots + U_n^{j_n+1}$$ 
and the following properties hold for a suitable constant $\Lambda>0$ only depending on $C, \Delta, D_k$ and $m$:
 \begin{enumerate}
  \item $\spt T_n^i$ and $\spt U_n^i$ are compact whenever $i\leq j_n$; furthermore, for each $i$ the sequence $(\spt T_n^i \cup \spt U_n^i)$, where $n$ is such that $j_n\geq i$, is uniformly compact;
  \item $\bdry T_n^2 = \ldots = \bdry T_n^{j_n+1} = 0$, and $\bdry T_n^1 = 0$
in case $m\ge 2$;
$$\sum_{i=1}^{j_n+1} \mass(T_n^i) < \Lambda;$$ 
if $m=1$ then $U_n^1 = \ldots = U_n^{j_n+1} = 0$, and if $m\geq 2$ then 
$$\sum_{i=1}^{j_n+1} \mass(U_n^i)+\mass(\bdry U_n^i) < \Lambda;$$ 
  \item for $1 \le L \le j_n-1$, the cycle $T_n^{L+1}+\ldots+T_n^{j_n+1}$
bounds an element of $\AKic_{m+1}(X_n)$ with mass less than $\Lambda\delta_L$, and $\mass(U_n^{L+1})+\ldots+\mass(U_n^{j_n+1})< \Lambda\sum_{i=L}^\infty \delta_i$.
 \end{enumerate}
\ep

For the proof of \thmref{thm:main} we further recall the definition of 
flat norm $\flatnorm(T)$ of an integral current $T\in\AKic_m(Z)$:
\begin{equation*}
 \flatnorm(T):= \inf\left\{\mass(U)+\mass(S) : T=U+\bdry S, U\in\AKic_m(Z), S\in\AKic_{m+1}(Z)\right\}.
\end{equation*}
A sequence $(T_n)$ in $\AKic_m(Z)$ converges in the flat topology
to a current $T \in \AKic_m(Z)$ if $\flatnorm(T - T_n) \to 0$.

\begin{proof}[Proof of \thmref{thm:main}]
 Choose numbers $0 < R_1< R_2 < \ldots \to\infty$ such that, after passing to a subsequence, we have $T_n\rstr\clB(x_n,R_r)\in\AKic_m(X_n)$ and
 \begin{equation*}
   \sup_n[\mass(T_n\rstr\clB(x_n,R_r)) + \mass(\bdry(T_n\rstr\clB(x_n,R_r)))] <\infty
 \end{equation*}
for $r\in\N$. Existence of such a sequence $(R_r)$ follows from \thmref{thm:slice-mass-bound} together with Fatou's Lemma, and the remark after \thmref{thm:bdry-rect}.
 Set $R_0:= 0$, and define $A_{r, n}:= \clB(x_n,R_r)\setminus \clB(x_n,R_{r-1})$ and $$T_{r,n}:= T_n\rstr A_{r,n}$$
for $r,n\in\N$; clearly $T_{r,n}\in\AKic_m(X_n)$ and
\begin{equation*}
   \sup_n[\mass(T_{r,n})+ \mass(\bdry T_{r,n})] <\infty.
 \end{equation*}
Let $T_{r,n}^1, \dots, T_{r,n}^{j_n+1},U_{r,n}^1, \dots, U_{r,n}^{j_n+1} \in\AKic_m(X_n)$ be currents as in \propref{prop:good-decomp-summary} for $T_{r,n}$ and $R_r$. For $n,s\in\N$, define closed sets
\begin{equation*}
 B_n^s:= \bigcup_{r=1}^s\bigcup_{i=1}^{\min\{s, j_n\}}(\{x_n\} \cup \spt T_{r,n}^i \cup \spt U_{r,n}^i)
\end{equation*}
and note that $B_n^1 \subset B_n^2 \subset \ldots \subset X_n$. According to part~(i) of \propref{prop:good-decomp-summary}, for each $s$, the sequence $(B_n^s)$ is uniformly compact. By~\cite[Proposition~5.2]{Wenger-cpt}, after passage to a subsequence, there exist isometric embeddings $\varphi_n\colon X_n\hookrightarrow Z$ and compact subsets $Y^1\subset Y^2\subset\ldots\subset Z$, for some complete metric space $Z$, such that 
\begin{equation*}
 \varphi_n(B_n^s)\subset Y^s
\end{equation*}
for all $n$ and $s$. Since $\varphi_n(x_n) \in Y^1$ for all $n$, we may arrange, by passing to a further subsequence, that $\varphi_n(x_n)$ converges to some $z_0\in Y^1$. 
Clearly, $\varphi_{n\#} U_{r,n}^i$ and $\varphi_{n\#} T_{r,n}^i$ are supported in $Y^s$ whenever $i \leq \min\{s,j_n\}$ and $r \leq s$. Moreover, for fixed $r$ and $i$, it follows from part~(ii) of \propref{prop:good-decomp-summary} that $\mass(\varphi_{n\#} T_{r,n}^i)$ and $\mass(\varphi_{n\#} U_{r,n}^i) + \mass(\bdry(\varphi_{n\#} U_{r,n}^i))$ are uniformly bounded and $\bdry(\varphi_{n\#} T_{r,n}^i)$ is either zero or, in case $m = 1$ and $i = 1$, equal to $\varphi_{n\#}(\bdry T_{r,n})$. We may therefore assume by the compactness and closure theorems in~\cite{AK-c}, after passing to a subsequence, that for every $r$ and $i$ there exist $T_r^i, U_r^i\in\AKic_m(Z)$ such that $$\varphi_{n\#} T_{r,n}^i \to T_r^i,\quad \varphi_{n\#} U_{r,n}^i\to U_r^i$$ weakly as $n\to\infty$. According to~\cite{Wenger-flatconv}, replacing $Z$ by $l^\infty(Z)$ if necessary, we may as well assume that the convergence is with respect to the flat topology.
Due to the lower semicontinuity of mass and assertion~(ii) of \propref{prop:good-decomp-summary}, we obtain that
\begin{equation*}
 \sum_{i=1}^\infty [\mass(T_r^i) + \mass(U_r^i) + \mass(\bdry U_r^i)]<\infty
\end{equation*}
and hence $\bar T_r:= \sum_{i=1}^\infty (T_r^i + U_r^i) \in \AKic_m(Z)$. Using part~(iii) of \propref{prop:good-decomp-summary}, one shows as in the last part of the proof of~\cite[Theorem~1.2]{Wenger-cpt} that, for every $r$,
\begin{equation}\label{eq:flatnorm-T_r}
  \flatnorm(\bar T_r-\varphi_{n\#}T_{r, n}) \to 0
\end{equation}
as $n\to\infty$. In particular, it follows that $\spt \bar T_r \subset \left\{z \in Z: R_{r-1} \le d(z_0,z) \le R_r\right\}$.

Now we view $T_{n,r}$ and $\bar T_r$ as elements of $\LWic{X_n}{m}$ and $\LWic{Z}{m}$, respectively. We define a function $T$ on $\Lipbs(Z)\times[\Liploc(Z)]^m$ by
\begin{equation*}
 T(f, \pi) := \sum_{r=1}^\infty \bar T_r(f, \pi),
\end{equation*}
where all but finitely many summands are zero because $\spt \bar T_r \cap \spt f = \emptyset$ for sufficiently large $r$. It is easily checked that $T\in\LWic{Z}{m}$. To show that $\varphi_{n\#}T_n \to T$ in the local flat topology, let $B \subset Z$ be a bounded closed set, and choose $s \in \N$ so that $\bar T^s := \bar T_1 + \ldots + \bar T_s$ satisfies $\|T - \bar T^s\|(B) = 0$, and also
$\|\varphi_{n\#}(T_n - T_n \rstr \clB(x_n,R_s))\|(B) = 0$.
It follows from \eqref{eq:flatnorm-T_r} that there exist $U_n\in\LWic{Z}{m}$ and $S_n\in\LWic{Z}{m+1}$ such that
\begin{equation*}
\bar T^s - \varphi_{n\#}(T_n \rstr \clB(x_n,R_s)) = U_n + \bdry S_n
\end{equation*}
and $\mass(U_n) + \mass(S_n)\to 0$. Now
\begin{equation*}
\|T - \varphi_{n\#}T_n - \bdry S_n\|(B) = 
\|\bar T^s - \varphi_{n\#}(T_n \rstr \clB(x_n,R_s)) - \bdry S_n\|(B),
\end{equation*}
hence $(\|T - \varphi_{n\#}T_n - \bdry S_n\| + \|S_n\|)(B)
= (\|U_n\| + \|S_n\|)(B) \to 0$.
This concludes the proof.
\end{proof}

\subsection{Local filling convergence}\label{subsection:fill-conv}

We now justify the remark after the statement of~\thmref{thm:main}. Replacing $Z$ by $l^\infty(Z)$ if necessary, we may assume that $Z$ admits isoperimetric inequalities of Euclidean type for $\AKic_k(Z)$, $k=1,\dots,m$. A consequence is the following useful fact, variations of which play a crucial role in the arguments of~\cite{Wenger-isopineq}, \cite{Wenger-flatconv}, \cite{Wenger-cpt}.

\bl \label{lem:spt-control}
For $k = 1,\dots,m+1$, there are constants $c_k$ such that, whenever $S \in \AKic_k(Z)$ and $\mass(S) < \delta^k$ for some $\delta > 0$, there exists $S' \in \AKic_k(Z)$ with $\bdry S' = \bdry S$, $\mass(S') < \delta^k$, and $d(x,\spt(\bdry S')) < c_k \delta$ for all $x \in \spt S'$.
\el

\begin{proof}
For $k \ge 2$, see~\cite[Lemma~3.4]{Wenger-isopineq}. For $k = 1$, a part of the argument is still valid. Given $S \in \AKic_1(Z)$ with $\mass(S) < \delta$ and a constant $Q > 1$, one gets a current $S' \in \AKic_1(Z)$ with $\bdry S' = \bdry S$ and $\mass(S') < \delta$ that is quasi-minimizing in the following sense: If $x \in \spt S'$, $0 < r < d(x,\spt(\bdry S'))$, and $S' \rstr \clB(x,r) \in \AKic_1(Z)$, then
\[
\mass(S' \rstr \clB(x,r)) \le Q \mass(Y)
\]
for every $Y \in \AKic_1(Z)$ with $\bdry Y = \bdry(S' \rstr \clB(x,r))$. Since $x \in \spt S'$, this shows in particular that the slice $\langle S',\varrho,r\rangle = \bdry(S' \rstr \clB(x,r)) \in \AKic_0(Z)$ with respect to the distance function $\varrho = d(x,\cdot)$ is non-zero, so that $\mass(\langle S',\varrho,r\rangle) \ge 2$. Integration from $0$ to $d(x,\spt(\bdry S'))$ gives $2 d(x,\spt(\bdry S')) \le \mass(S')$, hence $d(x,\spt(\bdry S')) < \delta/2$.
\end{proof} 

Suppose now that $(T_j)$ is sequence in $\LWic{Z}{m}$ that converges in the local flat topology to $0$, and suppose that for every bounded set $B \subset Z$, $\spt(\bdry T_j) \cap B = \emptyset$ for all but finitely many $j$. We want to show that then $T_j \to 0$ in the following sense: For every bounded closed set $B \subset Z$ there is a sequence $(S'_j)$ in $\LWic{Z}{m+1}$ such that $\spt(T_j - \bdry S'_j) \cap B = \emptyset$ for all but finitely many $j$, and $\|S'_j\|(B) \to 0$. This is an immediate consequence of the next result. We denote by $U_r(A)$ the open $r$-neighborhood of a set $A \subset Z$.

\bp
There is a constant $c > 0$ such that the following holds. 
Suppose $T \in \LWic{Z}{m}$, $B \subset Z$ is a
bounded closed set, $\delta > 0$, and $S \in \LWic{Z}{m+1}$ satisfies
$\|T - \bdry S\|(B) < \delta^m$ and $\|S\|(B) < \delta^{m+1}$.
Then there exists $S' \in \LWic{Z}{m+1}$ such that 
\[
\|T - \bdry S'\|(B) < \delta^m, \quad \mass(S') < c \delta^{m+1},
\] 
$\spt(T - \bdry S') \subset U_{c\delta}(\spt(\bdry T) \cup (Z \setminus B))$ 
and $\spt S' \subset U_{c\delta}(\spt T \cup (Z \setminus B))$.
\ep

\begin{proof}
Assume $B \ne \emptyset$. Put $R := T - \bdry S$, and fix $s > 0$ 
such that $\|R\|(U_s(B)) < \delta^m$ and $\|S\|(U_s(B)) < \delta^{m+1}$. 
Let $\varrho$ be the distance function to $B$.
There is an $r \in (0,s)$ such that, for $B_r := \{\varrho \le r\}$,
we have $T \rstr B_r, (\bdry S) \rstr B_r \in \AKic_m(Z)$ and 
$S \rstr B_r \in \AKic_{m+1}(Z)$. Then $R \rstr B_r \in \AKic_m(Z)$ and 
$\mass(R \rstr B_r) < \delta^m$. By \lemref{lem:spt-control},
there exists $R' \in \AKic_m(Z)$ such that
$\bdry R' = \bdry(R \rstr B_r)$, $\mass(R') < \delta^m$,
and $d(x,\spt(\bdry R')) < c_m\delta$ for all $x \in \spt R'$. 
Note that $\spt(\bdry R') \subset \spt(\bdry T) \cup (Z \setminus B)$. 
Since $R \rstr B_r - R'$ is a cycle with mass $< 2\delta^m$, 
the isoperimetric inequality of Euclidean type provides
a current $Q \in \AKic_{m+1}(Z)$ with $\bdry Q = R \rstr B_r - R'$ and
$\mass(Q) < 2^{(m+1)/m}D_m \delta^{m+1}$, for some constant $D_m$. 
Then $S \rstr B_r + Q \in \AKic_{m+1}(Z)$,
and $\mass(S \rstr B_r + Q) < (c'\delta)^{m+1}$ for some constant $c'$.
Using the above lemma again, we find $S' \in \AKic_{m+1}(Z)$ such that 
$\bdry S' = \bdry(S \rstr B_r + Q)$, $\mass(S') < (c'\delta)^{m+1}$,
and $d(x,\spt(\bdry S')) < c_{m+1}c'\delta$ for all $x \in \spt S'$. 
Note that $\bdry S' = \bdry(S \rstr B_r) + R \rstr B_r - R' = 
\langle S,\varrho,r\rangle + T \rstr B_r - R'$, so 
$T - \bdry S' = R' + T \rstr (Z \setminus B_r) - \langle S,\varrho,r\rangle$. 
Now $\|T - \bdry S'\|(B) = \|R'\|(B) < \delta^m$, 
and the result follows.
\end{proof}

\subsection{Uniqueness}\label{subsection:uniqueness}

We proceed to the discussion of~\propref{prop:intro-uniqueness}. We use a similar argument as in~\cite[Theorem~6.1]{Wenger-cpt}.

\begin{proof}[Proof of \propref{prop:intro-uniqueness}]
For each $n$, define the metric space $Z_n$ by gluing $Z$ and $Z'$ along $\varphi_n(X_n)$ and $\varphi'_n(X_n)$. Denote by $\varrho_n: Z\hookrightarrow Z_n$ and $\varrho'_n:Z'\hookrightarrow Z_n$ the natural isometric inclusions. Note that $$\varrho_n\circ\varphi_n = \varrho'_n\circ\varphi'_n$$ for all $n$.
Put $A := \{z_0\} \cup \spt T$ and $A' := \{z'_0\} \cup \spt T'$. Choose compact sets $C_1\subset C_2\subset\ldots\subset A$ and $C'_1\subset C'_2\subset\ldots\subset A'$  with $z_0 \in C_1$ and $z'_0 \in C'_1$ such that $\|T\|(A \setminus \bigcup C_i) = 0$ and $\|T'\|(A \setminus \bigcup C'_i) = 0$. Clearly, $\bigcup C_i$ is dense in $A$ and $\bigcup C'_i$ is dense in $A'$. Define
 \begin{equation*}
 B_n^i:=\varrho_n(C_i)\cup\varrho'_n(C'_i)
\end{equation*}
and note that $B_n^1\subset B_n^2\subset\ldots\subset Z_n$. Since $\varphi_n(x_n) \to z_0$ and $\varphi'_n(x_n) \to z'_0$, it follows that $d_{Z_n}(\varrho_n(z_0),\varrho'_n(z'_0)) \to 0$ as $n \to \infty$. For fixed $i$, the sequence $(B_n^i)$ is uniformly compact. By~\cite[Proposition 5.2]{Wenger-cpt}, we may assume,
after passing to a suitable subsequence, that there exist a complete metric space $Z''$, isometric embeddings $\sigma_n\colon Z_n\hookrightarrow Z''$, and compact subsets $Y^1\subset Y^2\subset\ldots\subset Z''$ such that
\begin{equation*}
 \sigma_n(B_n^i)\subset Y^i
\end{equation*}
for all $n$ and $i$. Consider the isometric embeddings $\tau_n:=\sigma_n\circ\varrho_n \colon Z \hookrightarrow Z''$ and $\tau'_n:= \sigma_n\circ\varrho'_n \colon Z' \hookrightarrow Z''$. Since $\tau_n(C_i)\subset Y^i$, we may assume, after passing to a subsequence, that $\tau_n\on{A}$ converges pointwise to an isometric embedding $\tau\colon A\hookrightarrow Z''$, uniformly on each $C_i$. Analogously, we may assume that $\tau'_n\on{A'}$ converges pointwise to an isometric embedding $\tau'\colon A'\hookrightarrow Z''$, uniformly on each $C'_i$. Since $d_{Z''}(\tau_n(z_0),\tau'_n(z'_0)) = d_{Z_n}(\varrho_n(z_0),\varrho'_n(z'_0)) \to 0$, we have 
\[
\tau(z_0) = \tau'(z'_0).
\]
It is not difficult to show that $\tau_{n\#}T\to \tau_\#T$ and $\tau'_{n\#}T'\to \tau'_\#T'$ weakly in $Z''$. We claim that also $\tau_{n\#}T - \tau'_{n\#}T' \to 0$ weakly. Then it follows that 
\begin{equation*}
 \tau_\#T - \tau'_\#T' = (\tau_\#T - \tau_{n\#}T) + (\tau_{n\#}T - \tau'_{n\#}T') + (\tau'_{n\#}T' - \tau'_\#T') \to 0
\end{equation*}
and thus $\tau_\# T = \tau'_\#T'$. Consequently, $\tau(\spt T) = \spt(\tau_\#T) = \spt(\tau'_\#T') = \tau'(\spt T')$, and $\psi:= \tau'^{-1}\circ\tau\colon (A,z_0) \to (A',z'_0)$ is a pointed isometry with $\psi_\#T=T'$.

To prove $\tau_{n\#}T - \tau'_{n\#}T' \to 0$, let first $B'' \subset Z''$ be a bounded closed set, and choose a bounded closed set $B \subset Z$ with $\tau_n^{-1}(B'') \subset B$ for all $n$; note that $\tau_n(z_0) \to \tau(z_0)$. Since $\varphi_{n\#}T_n \to T$ in the local flat topology, there is a sequence $(S_n)$ in $\LWic{Z}{m+1}$ such that 
$(\|T - \varphi_{n\#}T_n - \bdry S_n\| + \|S_n\|)(B) \to 0$, hence
\[
(\|\tau_{n\#}(T - \varphi_{n\#}T_n) - \bdry(\tau_{n\#} S_n)\| + 
\|\tau_{n\#} S_n\|)(B'') 
\le (\|T - \varphi_{n\#}T_n - \bdry S_n\| + \|S_n\|)(\tau_n^{-1}(B'')) \to 0.
\] 
This shows that $\tau_{n\#}(T - \varphi_{n\#}T_n) \to 0$ in the local flat topology of $\LWic{Z''}{m+1}$ and thus weakly. Analogously, $\tau'_{n\#}(T' - \varphi'_{n\#}T_n) \to 0$ weakly. Since $\tau_n \circ \varphi_n = \tau'_n \circ \varphi'_n$, we have
\[
\tau_{n\#}T - \tau'_{n\#}T' = \tau_{n\#}(T - \varphi_{n\#}T_n)   
+ \tau'_{n\#}(\varphi'_{n\#}T_n - T'),
\]
and the claim follows.
\end{proof}

\subsection{Ultralimits and Gromov--Hausdorff limits}\label{subsection:ultralimit}

It remains to prove \propref{prop:ultralimit}.
For the definitions of ultralimits and Gromov--Hausdorff limits of sequences of pointed metric spaces we refer to~\cite[Ch.~I.5]{Bridson-Haefliger} and~\cite[\S 8.1]{BBI}, respectively. 

\begin{proof}[Proof of \propref{prop:ultralimit}]
For every $z \in \spt T$ we choose a sequence $(y_n(z))$ with $y_n(z) \in \spt T_n$ such that $\varphi_n(y_n(z)) \to z$. This is clearly possible since $\varphi_{n\#}T_n \to T$ weakly and $\spt(\varphi_{n\#}T_n) \subset \varphi_n(\spt T_n)$. We have
\begin{equation}\label{eq:xy}
d_{X_n}(x_n,y_n(z)) = d_Z(\varphi_n(x_n),\varphi_n(y_n(z))) \to d_Z(z_0,z);
\end{equation}
furthermore, if $y'_n \in \spt T_n$ and $\varphi_n(y'_n) \to z' \in \spt T$, then 
\begin{equation}\label{eq:yy'}
d_{X_n}(y_n(z),y'_n) = d_Z(\varphi_n(y_n(z)),\varphi_n(y'_n)) \to d_Z(z,z'). 
\end{equation}
It follows that there is a well-defined isometric embedding $\psi \colon \{z_0\} \cup \spt T \to (X_\omega,x_\omega)$ that maps $z$ to the equivalence class $[(y_n(z))]$ of $(y_n(z))$ and $z_0$ to $[(x_n)] = x_\omega$. This proves~(i).

For part~(ii), since $\spt T$ is separable and $Y$ is proper, it suffices to show that for every finite set $F \subset \spt T$ there is an isometric embedding $f \colon \{z_0\} \cup F \to Y$ that maps $z_0$ to $y_0$. For every $z \in F$, choose a sequence $(y_n(z))$ as above, and let $E_n := \{x_n\} \cup \left\{y_n(z): z \in F\right\}$. Due to~\eqref{eq:xy}, there is an $r > 0$ such that $E_n \subset \clB(x_n,r)$ for every $n$. Hence, by the definition of pointed Gromov--Hausdorff convergence, there are maps $f_n \colon E_n \to Y$ such that $f_n(x_n) = y_0$ and 
\[
\max_{u,v \in E_n}|d_Y(f_n(u),f_n(v)) - d_{X_n}(u,v)| \to 0
\]
as $n \to \infty$. Since $Y$ is proper, we may assume that $f_n(y_n(z))$ converges to some $\bar y(z) \in Y$, for every $z \in F$. Then $d_{X_n}(x_n,y_n(z)) \to d_Y(y_0,\bar y(z))$ and $d_{X_n}(y_n(z),y_n(z')) \to d_Y(\bar y(z),\bar y(z'))$ for all $z,z' \in F$. Thus, in view of~\eqref{eq:xy} and~\eqref{eq:yy'}, we get an isometric embedding $f \colon \{z_0\} \cup F \to Y$ such that $f(z_0) = y_0$ and $f(z) = \bar y(z)$ for $z \in F$.
\end{proof}

Regarding the second part of \propref{prop:ultralimit}, note also that if a sequence of proper metric spaces $(X_n,x_n)$ converges to a complete metric space $(Y,y_0)$ in the pointed Gromov-Hausdorff sense, then clearly every bounded subset of $Y$ is totally bounded and hence $Y$ is proper.

\end{document}